\documentclass{gtart}

\def\ifplaintex{\expandafter\ifx\csname documentclass\endcsname\relax}

\expandafter\ifx\csname epsfbox\endcsname\relax\input epsf\fi

\ifplaintex 
\else
\fi

\def\gt{{\mathsurround=0pt\it $\cal G\mskip-2mu$eometry \&\ 
$\cal T\!\!$opology}}        

\def\gtp{{\mathsurround=0pt\it $\cal G\mskip-2mu$eometry \&\ 
$\cal T\!\!$opology $\cal P\!$ublications}}  

\def\lognumber#1{\def\thelognumber{#1}}
\def\volumenumber#1{\def\thevolumenumber{#1}}
\def\papernumber#1{\def\thepapernumber{#1}}
\def\volumeyear#1{\def\thevolumeyear{#1}}

\def\pagenumbers#1#2{\def\startpage{#1}\def\finishpage{#2}}
\def\published#1{\def\publishdate{#1}}
\def\proposed#1{\def\theproposer{#1}}
\def\seconded#1{\def\theseconders{#1}}
\def\received#1{\def\receiveddate{#1}}
\def\revised#1{\def\reviseddate{#1}}
\def\accepted#1{\def\accepteddate{#1}}

\long\def\asciiabstract#1{\long\def\theasciiabstract{#1}}
\def\asciikeywords#1{\def\theasciikeywords{#1}}

\let\\\par\let\thevolumenumber\relax\let\thepapernumber\relax
\let\thevolumeyear\relax\let\thesamplenumber\relax\let\startpage\relax
\let\finishpage\relax\let\publishdate\relax\let\receiveddate\relax
\let\reviseddate\relax\let\accepteddate\relax\let\theasciititle\relax
\let\theasciiauthors\relax
\let\theasciiabstract\relax\let\theasciikeywords\relax
\let\theasciiemail\relax\let\theshortauthors\relax\let\theshorttitle\relax

\long\def\maketitlep{   

\count0=\startpage

\gt\hfill      
\hbox to 77pt{\vbox to 0pt{\vglue -15pt\epsfbox{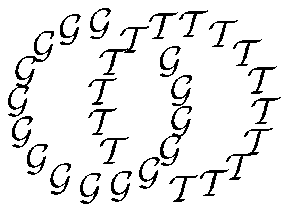}\vss}\hss}
\break
{\small\ifx\thesamplenumber\relax 
Volume \else Sample
\fi\thevolumenumber\ (\thevolumeyear)
\startpage--\finishpage\nl
Published: \publishdate}
\vglue 0.5truein plus 0.4fil minus 0.1truein

{\parskip=0pt\leftskip 0pt plus 1fil\def\\{\par\smallskip}{\ifplaintex\large
\else\Large\fi\bf\thetitle}\par\medskip}   

\vglue 0pt plus 0.1fil 

{\parskip=0pt\leftskip 0pt plus 1fil\def\\{\par}{\sc\theauthors}
\par\medskip}

\vglue 0pt plus 0.1fil 

{\small\parskip=0pt\let\newline\\
{\leftskip 0pt plus 1fil\def\\{\par}{\sl\theaddress}\par}
\expandafter\ifx\theemail\relax    
\relax\else\vglue 5pt plus 0.02fil minus 2pt\def\\{\stdspace{\rm 
and}\stdspace} 
\cl{Email:\stdspace\tt\theemail}\fi
\ifx\theurl\relax                  
\relax\else\vglue 5pt plus 0.02fil minus 2pt\def\\{\stdspace{\rm 
and}\stdspace}
\cl{URL:\stdspace\tt\theurl}\fi\par}

\vglue 7pt plus 0.3fil minus 3pt

{\bf Abstract}
\vglue 5pt plus 0.1fil minus 2pt

\theabstract

\vglue 7pt plus 0.3fil minus 3pt

{\bf AMS Classification numbers}\quad Primary:\quad \theprimaryclass

Secondary:\quad \thesecondaryclass

\vglue 5pt plus 0.3fil minus 2pt

{\bf Keywords:}\quad \thekeywords

\vglue 10pt plus 0.5fil minus 5pt

{\small  Proposed: \theproposer\hfill Received: \receiveddate\nl
Seconded: \theseconders\hfill 
\ifx\reviseddate\relax                         
Accepted: \accepteddate                        
\else
Revised: \reviseddate                          
\fi}
\eject
}       

\let\maketitlepage\maketitlep
\let\maketitle\maketitlepage

\font\phead=cmsl9 scaled 950
\font=cmsl9 scaled 1050
\font\pnum=cmbx10 scaled 913
\font\lnum=cmbx10 
\font\pfoot=cmsl9 scaled 950
\font=cmsl9 scaled 1050
\ifplaintex
\headline{\vbox to 0pt{\vskip -4.5mm\line{\small\phead\ifnum
\count0=\startpage ISSN 1364-0380 (on line)
1465-3060 (printed) \hfill {\pnum\folio}\else\ifodd\count0\def\\{ }%
\ifx\theshorttitle\relax\thetitle\else\theshorttitle\fi\hfill{\pnum\folio}
\else\def\\{ and }{\pnum\folio}\hfill\ifx\theshortauthors\relax\theauthors
\else\theshortauthors\fi\fi\fi}\vss}}
\footline{\vbox to 0pt{\vglue 0mm\line{\small\pfoot\ifnum\count0=\startpage
\copyright\ \gtp\hfill\else
\gt, Volume \thevolumenumber\ (\thevolumeyear)\hfill\fi}\vss
}}
\else
\makeatletter
\def\@oddhead{{\small\lhead\ifnum\count0=\startpage ISSN 1364-0380 (on line)
1465-3060 (printed) \hfill {\lnum\number\count0}\else\ifodd\count0
\def\\{ }\ifx\theshorttitle\relax \thetitle \else\theshorttitle\fi\hfill
{\lnum\number\count0}\else\def\\{ and }{\lnum\number\count0}
\hfill\ifx\theshortauthors\relax 
\theauthors\else\theshortauthors\fi\fi\fi}}\def\@evenhead{@oddhead}
\def\@oddfoot{\small\lfoot\ifnum\count0=\startpage\copyright\ \gtp\hfill\else
\gt, Volume \thevolumenumber\ (\thevolumeyear)\hfill\fi}
\def\@evenfoot{@oddfoot}
\makeatother
\fi

\newwrite\gtoutfile
\long\gdef\makeheadfile{  
{\def\\{, }\def\s{ }
\immediate\openout\gtoutfile head.xxx
\immediate\write\gtoutfile{To: math@arxiv.org}
\immediate\write\gtoutfile{Subject: put OR rep NNNNN:pppp}
\immediate\write\gtoutfile{--text follows this line--}
\immediate\write\gtoutfile{Proxy-for: \ifx\theasciiauthors\relax
\theauthors\else\theasciiauthors\fi\s<\ifx\theasciiemail\relax\theemail\else\theasciiemail\fi>}
\immediate\write\gtoutfile{\noexpand\\}
\immediate\write\gtoutfile{Authors: \ifx\theasciiauthors\relax
\theauthors\else\theasciiauthors\fi}
{\def\\{ }\immediate\write\gtoutfile{Title: \ifx\theasciititle\relax
\thetitle\else\theasciititle\fi}}
\immediate\write\gtoutfile{Subj-class: GT or GR or SG or ...}
\immediate\write\gtoutfile{MSC-class: \theprimaryclass\ifx\thesecondaryclass\relax\else, \thesecondaryclass\fi}
\immediate\write\gtoutfile{Journal-ref: Geom. Topol. \thevolumenumber\s
(\thevolumeyear) \startpage-\finishpage}
\immediate\write\gtoutfile{Comments: Published in Geometry and Topology at}
\immediate\write\gtoutfile{    http://www.maths.warwick.ac.uk/gt/GTVol\thevolumenumber/paper\thepapernumber.abs.html}
\immediate\write\gtoutfile{\noexpand\\}
\immediate\write\gtoutfile{}
\ifx\theasciiabstract\relax
\immediate\write\gtoutfile{\theabstract}\else
\immediate\write\gtoutfile{\theasciiabstract}\fi
\immediate\write\gtoutfile{}
\immediate\write\gtoutfile{\noexpand\\}
\immediate\write\gtoutfile{}
\immediate\closeout\gtoutfile}}  

\def\maketitlepage{\maketitlep\makeheadfile}
\let\maketitle\maketitlepage

\lognumber{190}
\volumenumber{6}\papernumber{8}\volumeyear{2002}
\pagenumbers{219}{267}
\received{21 June 2001}
\revised{21 March 2002}
\published{3 May 2002}
\accepted{3 May 2002}

\proposed{Walter Neumann}
\seconded{Cameron Gordon, Wolfgang Metzler}

\usepackage{amsmath,amssymb,amsxtra}

\def\S{Section }

\makeatletter
\newtheoremstyle{noname}{14pt plus6.3pt minus6.3pt}{7.4pt plus3pt minus3pt}%
{\rm}{}{\bf}{}{0.75em}{\thmnumber{#2}\thmnote{\sl{\stdspace#3}}}
  \def\tagform@#1{\maketag@@@{%
   \textbf{(\ignorespaces#1\unskip\@@italiccorr)}}}%
   \renewcommand{\eqref}[1]{\textup{\maketag@@@{(\ignorespaces%
	{\ref{#1}}\unskip\@@italiccorr)}}}
\makeatother

\theoremstyle{plain}
\newtheorem{theorem}[equation]{Theorem}
\newtheorem{lemma}[equation]{Lemma}
\newtheorem{corollary}[equation]{Corollary}
\newtheorem{proposition}[equation]{Proposition}
\newtheorem{mainthm}[equation]{Theorem}
\newtheorem{maincor}[equation]{Corollary}

\theoremstyle{definition}
\newtheorem{definition}[equation]{Definition}
\newtheorem{remark}[equation]{Remark}
\newtheorem*{ack}{Acknowledgements}

\newtheorem{question}[equation]{Question}
\newtheorem*{notation}{Notation}

\theoremstyle{noname}
\newtheorem{subsec}[equation]{}
\newtheorem{subexample}{Example}[equation]

\numberwithin{equation}{section}
\renewcommand{\theenumi}{\alph{enumi}}

\newenvironment{proof2}%
  {\begin{proof}}{\renewcommand{\qed}{}\end{proof}}

\DeclareFontFamily{OML}{rsfs}{\skewchar\font'177}
\DeclareFontShape{OML}{rsfs}{m}{n}{ <5> <6> rsfs5 <7> <8> <9> rsfs7
  <10> <10.95> <12> <14.4> <17.28> <20.74> <24.88> rsfs10 }{}
\DeclareMathAlphabet{\mathfs}{OML}{rsfs}{m}{n}

\newcommand{\As}{{\mathfs A}}
\newcommand{\Bs}{{\mathfs B}}
\newcommand{\Es}{{\mathfs E}}
\newcommand{\Z}{{\mathbb Z}}
\newcommand{\Q}{{\mathbb Q}}
\newcommand{\abs}[1]{\left\lvert#1\right\rvert}
\renewcommand{\leq}{\leqslant}
\renewcommand{\geq}{\geqslant}
\renewcommand{\emptyset}{\varnothing}
\DeclareMathOperator{\id}{id}
\DeclareMathOperator{\ad}{ad}
\DeclareMathOperator{\dist}{dist}

\input xy
\xyoption{all}

\begin{document}

\title{Deformation and rigidity of simplicial\\group actions on trees} 
\author{Max Forester}
\address{Mathematics Institute, University of Warwick\\Coventry, CV4 7AL,
UK} 
\email{forester@maths.warwick.ac.uk}
\begin{abstract}
We study a notion of deformation for simplicial trees with group actions
($G$--trees). Here $G$ is a fixed, arbitrary group. Two $G$--trees are
related by a deformation if there is a finite sequence of collapse and
expansion moves joining them. We show that this relation on the set of
$G$--trees has several characterizations, in terms of dynamics, coarse
geometry, and length functions. Next we study the deformation space of a
fixed $G$--tree $X$. We show that if $X$ is ``strongly slide-free'' then
it is the unique reduced tree in its deformation space. 

These methods allow us to extend the rigidity theorem of Bass and
Lubotzky to trees that are not locally finite. This yields a unique
factorization theorem for certain graphs of groups. We apply the theory
to generalized Baumslag--Solitar groups and show that many have canonical
decompositions. We also prove a quasi-isometric rigidity theorem for
strongly slide-free $G$--trees. \end{abstract}

\asciiabstract{
We study a notion of deformation for simplicial trees with group actions
(G-trees). Here G is a fixed, arbitrary group. Two G-trees are
related by a deformation if there is a finite sequence of collapse and
expansion moves joining them. We show that this relation on the set of
G-trees has several characterizations, in terms of dynamics, coarse
geometry, and length functions. Next we study the deformation space of a
fixed G-tree X. We show that if X is `strongly slide-free' then
it is the unique reduced tree in its deformation space. 

These methods allow us to extend the rigidity theorem of Bass and
Lubotzky to trees that are not locally finite. This yields a unique
factorization theorem for certain graphs of groups. We apply the theory
to generalized Baumslag-Solitar groups and show that many have canonical
decompositions. We also prove a quasi-isometric rigidity theorem for
strongly slide-free G-trees.} 

\asciikeywords{G-tree, graph of groups, folding, Baumslag-Solitar group,
quasi-isometry} 
\maketitlepage

\section{Introduction}

In this paper we study a notion of deformation for simplicial trees with
group actions ($G$--trees). Here $G$ is a fixed, arbitrary group, and we
consider actions that do not invert edges. Such actions correspond to
graph of groups decompositions of $G$, according to Bass--Serre
theory \cite{serre:trees}. 

Our notion of deformation is based on \emph{collapse moves} in graphs of
groups, in which an edge carrying an amalgamation of the form $A \ast_C
C$ is collapsed to a vertex with group $A$. This operation simplifies the
underlying graph without enlarging any vertex or edge
groups. Collapse moves can be defined and performed
directly on $G$--trees as well. An \emph{elementary deformation} is a finite
sequence of collapse moves and their inverses (called \emph{expansion
moves}). We are interested in knowing when two $G$--trees are related by
a deformation, and in what can be said about them if they are. These two
questions are addressed in the two main theorems of the paper. 

In order to discuss the first result we recall some definitions. If $X$
is a $G$--tree then an \emph{elliptic subgroup} is a subgroup $H
\subseteq G$ that fixes a vertex of $X$. Thus, the elliptic subgroups are
precisely the vertex stabilizers and their subgroups. The \emph{length
function} $\ell_X\co G \to \Z$ assigns to $\gamma\in G$ the minimum
displacement of any vertex under $\gamma$. Here we regard $G$--trees as
metric spaces by assigning length one to every edge. A
\emph{quasi-isometry} is a map which preserves the large scale geometry of
a metric space. The definition is fairly standard; see section
\ref{geomsec} for details. In this section we also define coarse
equivariance, a slight weakening of the property of
equivariance. Collecting together Theorem \ref{mainthm}, Corollary
\ref{mainthmcor}, and Theorem \ref{qicor}, we have the first main result. 

\begin{mainthm} \label{maindefthm}
Let $G$ be a group and let $X$ and $Y$ be cocompact $G$--trees. The
following conditions are equivalent. 
\begin{enumerate} 
\item \label{a1} $X$ and $Y$ are related by an elementary deformation. 
\item \label{a2} $X$ and $Y$ have the same elliptic subgroups. 
\item \label{a3} There is a coarsely equivariant quasi-isometry $\phi \co
X \to Y$. 
\end{enumerate}
If all vertex stabilizers are finitely generated then we may also
include: 
\begin{enumerate}
\setcounter{enumi}{3}
\item \label{a4} The length functions $\ell_X$ and $\ell_Y$ vanish on the
same elements of $G$. 
\end{enumerate}
\end{mainthm}

Note that condition (\ref{a2}) does not require $X$ and $Y$ to have the
same vertex stabilizers. In fact, trees related by a deformation can
easily have different vertex stabilizers.  Condition (\ref{a2}) arises
naturally in many situations, making the implication
(\ref{a2})$\Rightarrow$(\ref{a1}) quite useful. 

Regarding condition (\ref{a4}), it is well known that in most cases a
$G$--tree is determined by its length function \cite{cullermorgan}.  Our
result shows that if one has partial knowledge of $\ell_X$ (namely, its
vanishing set), then $X$ is partly determined, in an understandable way:
up to deformation or up to quasi-isometry.  The implication
(\ref{a3})$\Rightarrow$(\ref{a1}) is also interesting, as it transforms a
coarse, approximate relationship between $G$--trees into a precise one. 

One direct application of Theorem \ref{maindefthm} is a proof of the
conjecture of Herrlich stated in \cite{herrlich} (see Corollary
\ref{herrlichcor}). For other applications it is helpful to know more
about condition (\ref{a1}). This is the subject of the second result. We
fix a $G$--tree $X$ and consider the set of all trees related to $X$ by
elementary deformations. This set is called the \emph{deformation space}
of $X$. We say that a $G$--tree is \emph{reduced} if it admits no
collapse moves. Every cocompact $G$--tree can be made reduced by
performing collapse moves. 

A $G$--tree is \emph{strongly slide-free} if its stabilizers satisfy the
following condition: for any edges $e$ and $f$ with common initial vertex
$v$, if $G_e \subseteq G_f$ then there is an element $\gamma \in G$
fixing $v$ and taking $e$ to $f$. In terms of graphs of groups this means
that for every vertex group $A$, if $C$ and $C'$ are neighboring edge
groups then no conjugate (in $A$) of $C$ is contained in $C'$. Our second
main result shows that strongly slide-free trees are locally rigid, from
the point of view of deformations. (See Theorem \ref{localrigidity} for a
more comprehensive statement.) 

\begin{mainthm} \label{mainrigidthm} 
Let $X$ and $Y$ be cocompact $G$--trees that are related by an elementary
deformation. If $X$ is strongly slide-free and $Y$ is reduced, then there
is a unique $G$--isomorphism $X \to Y$.
\end{mainthm}

This means that if a deformation space contains a strongly slide-free
$G$--tree, then this tree is the unique reduced $G$--tree in the
space. Hence, all trees in the deformation space reduce (by collapse
moves) to the same $G$--tree. 
Combining Theorems \ref{maindefthm} and \ref{mainrigidthm} we obtain the
following rigidity theorem. 

\begin{maincor} \label{rigiditycor} 
Let $X$ and $Y$ be cocompact $G$--trees with the same elliptic
subgroups. If $X$ is strongly slide-free and $Y$ is reduced then there
is a unique isomorphism of $G$--trees $X \to Y$. 
\end{maincor}

This result provides an answer to the question of Bass and Lubotzky
that is raised in \cite{bass:rigidity}. It generalizes their Rigidity
Theorem, which applies only to locally finite trees. A statement that
more closely resembles theirs is given in Corollary \ref{rigidcor}. In
fact this latter result, expressed in graph of groups language, has an
interesting application. Recall that a group is \emph{unsplittable} if it
admits no nontrivial graph of groups decomposition. Then we have: 

\begin{maincor} \label{uniquefactorcor} 
Every group is the fundamental group of at most one strongly slide-free
graph of unsplittable finitely generated groups, with finite underlying
graph. 
\end{maincor}

The phrase ``at most one'' is meant up to graph of groups isomorphism in
the sense of \cite{bass:covering}. This result is somewhat analogous to
the classical theorem that states that every group has at most one free
product decomposition (up to rearrangement of factors) into groups that
are not themselves free products \cite[\S 35]{kurosh:groups}. 

Our results can be applied to the study of \emph{generalized
Baumslag--Solitar trees}, which are $G$--trees whose vertex and edge
stabilizers are all infinite cyclic. The groups $G$ that arise are called
\emph{generalized Baumslag--Solitar groups}. They include the classical
Baumslag--Solitar groups, torus knot groups, and finite index subgroups
of these groups. They have the virtue that their elliptic subgroups are
uniquely determined, independently of the tree (except in some degenerate
cases). Thus for any such group $G$, there is a single deformation space
of $G$--trees. If there is a strongly slide-free $G$--tree then it is
canonical, by Theorem \ref{mainrigidthm}. A general statement for
generalized Baumslag--Solitar trees is given in Corollary \ref{gbscor}. 

A geometric application of Theorem \ref{mainrigidthm} is obtained using
condition (\ref{a3}) of Theorem \ref{maindefthm}. It is a quasi-isometric
rigidity theorem for trees. This result complements the work of Mosher,
Sageev, and Whyte \cite{msw:announcement}, though there is no direct
connection between their work and ours; they study quasi-isometries of
groups, whereas here we work with a fixed group. 

\begin{maincor} \label{geomrigiditycor} 
Let $X$ and $Y$ be reduced cocompact $G$--trees with $X$ strongly
slide-free. Given a coarsely equivariant quasi-isometry $\phi \co
X \to Y$, there is a unique equivariant isometry $X \to Y$, and it
has finite distance from $\phi$. 
\end{maincor}

Thus, strongly slide-free cocompact $G$--trees are quasi-isometrically
rigid, in the equivariant sense. This means in particular that the local
geometry of $X$ is completely determined by its equivariant large scale
geometry. 

Finally we mention that Theorems \ref{maindefthm} and \ref{mainrigidthm}
can be used to obtain uniqueness results for various decompositions of
groups, such as one-ended decompositions of accessible groups and JSJ
decompositions.  These applications are described in detail in
\cite{forester:jsj}. 

\medskip

Our approach to proving Theorem \ref{maindefthm} is based on the method of
folding of $G$--trees. This technique was used by Chiswell in 
his thesis (see \cite{chiswell}), and has also been developed by
Stallings, Bestvina and Feighn, Dunwoody, and others
\cite{stallings:foldings,bestvina:accessibility,dunwoody:folding}. In
this paper, in order to understand the most general situation, we define
and analyze folds performed at infinity, or \emph{parabolic folds}. 

During such a move, various rays with a common end become identified.
The simplest example is given by the quotient map of a parabolic tree, in
which edges or vertices are identified if they have the same relative
distance from the fixed end. The result is a linear
tree. More generally, fewer identifications may be made, so that the
result is ``thinner'' parabolic tree, or similar operations may
be performed inside a larger $G$--tree.
Under suitable conditions, a morphism between $G$--trees can be factored
as a finite composition of folds, multi-folds, and parabolic folds. We
show that such a morphism can be constructed when two $G$--trees have the
same elliptic subgroups. Then the various types of folds are shown to be
elementary deformations, in this particular situation.

The proof of the Theorem \ref{mainrigidthm} relies on the notion of a
\emph{telescoping} of a $G$--tree. A telescoping is a structure on a
$G$--tree $X$ which remembers a second $G$--tree from which $X$
originated. We show that if one begins with a strongly slide-free
$G$--tree, then this structure is preserved by elementary deformations,
and so a ghost of the original tree is always present throughout any 
deformation. The theorem is proved by observing that any maximal sequence
of collapse moves will then recover the original tree.

\begin{ack} 
I would like to thank Koji Fujiwara and Peter Scott for helpful
conversations related to this work. I also thank David Epstein 
for his encouragement. This work was supported by EPSRC grant
GR/N20867. 
\end{ack}

\section{Basic properties of $G$--trees} 

A \emph{graph} $A = (V(A),E(A))$ is a pair of sets
together with a fixed point free involution $e \mapsto \overline{e}$ of
$E(A)$ and maps $\partial_0, \partial_1 \co E(A) \rightarrow V(A)$ such
that $\partial_i(\overline{e}) = \partial_{1-i}(e)$ for every $e\in
E(A)$. Elements of $V(A)$ are called \emph{vertices} and elements of
$E(A)$ are called \emph{edges}. The pair $\{e, \overline{e}\}$ is called
a \emph{geometric edge}. An edge $e$ for 
which $\partial_0 e = \partial_1 e$ is called a \emph{loop}. For each
vertex $v\in V(A)$ we set $E_0(v) = \{e\in E(A) \mid \partial_0 e = v \}$. 

A \emph{graph of groups} ${\bf A} = (A, \As, \alpha)$ consists of the
following data: a connected graph $A$, groups $\As_a$ ($a \in V(A)$) and
$\As_e =  \As_{\overline{e}}$ ($e \in E(A)$), and injective homomorphisms
$\alpha_e \co \As_e \rightarrow \As_{\partial_0 e}$ ($e \in E(A)$). Given
a vertex $a_0\in V(A)$ there is a \emph{fundamental group} $\pi_1({\bf
A},a_0)$ whose isomorphism type does not depend on $a_0$ (see
\cite[Chapter I, \S 5.1]{serre:trees}). 

A \emph{tree} is a connected graph with no circuits: if $(e_1, \ldots,
e_n)$ is a path with $e_{i+1} \not= \overline{e}_i$ for all $i$ (ie, a
path \emph{without reversals}), then $\partial_0 e_1 \not= \partial_1
e_n$. 

A \emph{ray} is a semi-infinite path without reversals. An \emph{end} of
a tree $X$ is an equivalence class of rays, where two rays are considered
equivalent if their intersection is a ray. The set of ends is called the
\emph{boundary} of the tree, denoted $\partial X$. A subtree
\emph{contains} an end if it contains a representative of that end. 

An automorphism of a tree is an \emph{inversion} if it maps $e$ to
$\overline{e}$ for some edge $e$. 

\begin{definition} 
Let $G$ be a group. A \emph{$G$--tree} is a tree $X$ together with an
action of $G$ on $X$ by automorphisms, none of which are
inversions. There is a well known correspondence between $G$--trees and
graphs of groups having fundamental group $G$. A good reference for this
material is \cite{bass:covering} (see also \cite{serre:trees}). 
\end{definition}

\begin{definition}
A $G$--tree is \emph{minimal} if it contains no proper invariant
subtree. Correspondingly, a graph of groups is minimal if there is no
proper subgraph $A' \subsetneq 
A$ which carries the fundamental group of ${\bf A}$. That is, the
injective homomorphism 
 $\pi_1(({\bf A}\vert_{A'}), a_0) \rightarrow \pi_1({\bf A},a_0)$ 
of fundamental groups induced by the inclusion of 
${\bf A} \vert_{A'}$ into ${\bf A}$ 
is not surjective for any $A' \subsetneq A$. 

If $A$ has finite diameter then this condition is equivalent to the
following property: for every vertex $v\in V(A)$ of valence one, the
image of the neighboring edge group is a proper subgroup of the vertex
group $\As_v$ (see \cite[7.12]{bass:covering}). 
\end{definition}

\begin{definition}\label{reduced}
A $G$--tree is \emph{reduced} if, whenever $G_e = G_{\partial_0 e}$ for
an edge $e$ of the tree, $\partial_0 e$ and $\partial_1 e$ are in the
same $G$--orbit. This occurs if and only if the tree admits no ``collapse
moves'' (see \ref{collapse} below). The corresponding notion for graphs
of groups is: for every edge $e\in E(A)$, if $\alpha_e(\As_e) =
\As_{\partial_0 e}$ then $e$ is a loop. 

This definition differs from the notion of ``reduced'' used in
\cite{bestvina:accessibility}. Note that if a $G$--tree (or a graph of
groups) is reduced then it is minimal. 

Every cocompact $G$--tree can be made reduced by performing collapse
moves until none are available. These moves are described in the next
section. The resulting $G$--tree will depend on the particular sequence
of collapse moves chosen. Similarly, graphs of groups with finite
underlying graphs can be made reduced. 
\end{definition}

\begin{definition} \label{ellipticsubgroup} 
Let $X$ be a $G$--tree. An element $\gamma \in G$ is \emph{elliptic} if it
has a fixed point, and \emph{hyperbolic} otherwise. We define the
\emph{length function} of $X$ by
\begin{equation*}
\ell_X(\gamma) = \min_{x\in V(X)} d(x,\gamma x). 
\end{equation*}
Thus, $\ell_X(\gamma) = 0$ if and only if $\gamma$ is elliptic. If
$\gamma$ is hyperbolic then the set
\[ L_{\gamma} \ = \ \{ x \in X \mid d(x,\gamma x) = \ell_X(\gamma) \} \]
is a $\gamma$--invariant linear subtree, called the \emph{axis} of
$\gamma$. The action of $\gamma$ on its axis is by a
translation of amplitude $\ell_X(\gamma)$.

A \emph{minimal subtree} is a nonempty $G$--invariant subtree which is
minimal. If $G$ contains a hyperbolic element then there is a unique
minimal subtree, equal to the union of the axes of hyperbolic elements
(see \cite[7.5]{bass:covering}). If $G$ contains no hyperbolic elements
then any global fixed point is a minimal subtree. It is possible for a
$G$--tree to have no minimal subtree. 

A $G$--tree is \emph{elliptic} if there is a global fixed point. It is
\emph{parabolic} if there is a fixed end, and some element of $G$ is
hyperbolic. It follows that the minimal subtree is also parabolic, and
has quotient graph equal to a closed circuit. A parabolic $G$--tree may
have two fixed ends. In this case the minimal subtree is a linear tree
acted on by translations. Otherwise the fixed end is unique. (Any
$G$--tree with three fixed ends is elliptic.) See \cite[\S
2]{cullermorgan} for a more complete discussion of these facts. 

\begin{remark} \label{parabolicsubgroup} 
In a parabolic $G$--tree
the set of elliptic elements is a subgroup, whereas this conclusion is
false for $G$--trees in general (cf Lemma \ref{hypsegment}). 
To see this, let $\varepsilon$ be a fixed end. Every elliptic element
fixes pointwise a ray tending to $\varepsilon$. Since any two such rays
have nonempty intersection, any two elliptic elements have a common fixed
point. Thus, their product is elliptic. 
\end{remark}

A subgroup $H$ of $G$ is
\emph{elliptic} if it fixes a vertex of $X$ (ie, if $X$ is an elliptic
$H$--tree). This property is stronger than requiring the elements of
$H$ to be elliptic, though in many cases these two properties coincide
(cf Proposition \ref{titslemma}). 
\end{definition}

\begin{notation}
If $\gamma \in G$ and $H\subseteq G$ is a subgroup, set $H^{\gamma} =
\gamma H \gamma^{-1}$. Note that this yields the identity
$(H^{\delta})^{\gamma} = H^{(\gamma \delta)}$. If $x$ is a vertex or edge
of a $G$--tree then the \emph{stabilizer} of $x$ is $G_x = \{\gamma\in
G \mid \gamma x = x\}$. All group actions are on the left, so that 
$G_{\gamma x} = (G_x)^{\gamma}$. 

If $x$ and $y$ are vertices, edges, or ends of a tree, let $[x,y]$ denote
the unique smallest subtree containing $x$ and $y$. It is an unoriented
segment (including edges and their inverses), possibly with zero or
infinite length. In the case where $y$ is an end, we will sometimes use
the more suggestive notation $[x,y)$. 
Two edges $e$ and $f$ are \emph{coherently oriented} if
they are members of an oriented path without reversals. This occurs if
and only if the segment $[\partial_0 e, \partial_0 f]$ contains exactly
one of $e$, $f$.
\end{notation}

\begin{proposition}[Tits {\cite[3.4]{tits:treeauto}}]
\label{titslemma} 
Let $X$ be a $G$--tree. If every element of $G$ has a fixed point in $X$
then either there is a global fixed point, or there is a unique end
$\varepsilon \in \partial X$ which is fixed by $G$.  In the latter case,
if $(x_1,x_2, x_3, \ldots)$ is a sequence of vertices or edges tending
monotonically to $\varepsilon$, then $G_{x_i} \subseteq G_{x_{i+1}}$ for
all $i$, with strict inclusion for infinitely many $i$, and $G =
\bigcup_{i \geq 0} G_{x_i}$. \endproof
\end{proposition}

\begin{lemma}\label{hypcriteria}
Let $X$ be a $G$--tree. An element $\gamma \in G$ is
hyperbolic under any one of the following conditions: 
\begin{enumerate}
\item \label{h1} $d(x, \gamma x)$ is odd for some vertex $x$; 
\item \label{h2} for some vertex $x$ not fixed by $\gamma$, the edges of
the path $[x, \gamma x]$ map injectively to the edges of $G
\backslash X$; 

\item \label{h3} for some edge $e$ not fixed by $\gamma$, the
edges $e$ and $\gamma e$ are coherently oriented. 
\end{enumerate}
\end{lemma}

\begin{proof}
If $\gamma$ is not hyperbolic then it fixes some vertex $v\in V(X)$. 
To rule out cases (\ref{h1}) and (\ref{h2}), suppose that $x$ is any
vertex not fixed by $\gamma$. The
subtree $T$ spanned by $v$, $x$, and $\gamma x$ is equal to $[v,x] \cup
[v,\gamma x]$, and $[v,x] \cap [v,\gamma x]$ meets $[x,\gamma x]$ in a
single vertex, $w$. This vertex is fixed by $\gamma$ and so $\gamma
([w,x])$ = $[w,\gamma x]$. As $[x,\gamma x] = [x,w] \cup [w,\gamma x]$,
this path has even length and contains a pair of edges which are related
by $\gamma$. These two edges have the same image in $G \backslash 
X$. This shows that neither (\ref{h1}) nor (\ref{h2}) can hold. 

To rule out case (\ref{h3}) let $e$ be the edge, and look at the subtree
spanned by $e$, $\gamma e$, and $v$. Then as $\gamma([v,e]) = [v,\gamma
e]$, the edges $e$ and $\gamma e$ are both oriented toward $v$, or both
oriented away from $v$. Since any oriented path traversing $[e, \gamma
e]$ travels first toward $v$ and then away from $v$, the two edges are
incoherently oriented, violating (\ref{h3}). 
\end{proof}

\begin{lemma}[Hyperbolic Segment Condition] \label{hypsegment}
Let $X$ be a $G$--tree and let $e$ and $f$ be edges such that $[e,f] =
[\partial_0 e, \partial_0 f]$. Given $\gamma_e \in G_{\partial_0 e} -
G_e$ and $\gamma_f \in G_{\partial_0 f} - G_f$, the product $\gamma_e
\gamma_f$ is hyperbolic and its axis contains $[e,f]$. 
\end{lemma}

\begin{proof}
The element $\gamma_e \gamma_f$ takes $\gamma_f^{-1}(e)$ to $\gamma_e(e)$
and one easily checks that these two edges are coherently oriented. Then
$\gamma_e \gamma_f$ is hyperbolic by Lemma
\ref{hypcriteria}(\ref{h3}). The subtree $\bigcup_{n \in \Z} (\gamma_e
\gamma_f)^n [\gamma_f^{-1}(e),e]$ is a linear $\gamma_e
\gamma_f$--invariant subtree, and must therefore be the axis. Now note that
$[e,f] \subseteq [\gamma_f^{-1}(e),e]$. 
\end{proof}

\section{Moves and factorizations}

In this section we discuss ways of modifying a $G$--tree to obtain
another. We also discuss how to factor the more complicated moves as
compositions of simpler moves, under certain conditions. The simplest of
these moves will be called \emph{elementary moves}. None of the moves
affect the group $G$. In particular, a sequence of moves relating two
graphs of groups will induce an isomorphism between their fundamental
groups. 

In all of the descriptions below, $X$ is a $G$--tree and ${\bf A} =
(A, \As, \alpha)$ is the corresponding graph of groups. 

\begin{remark} \label{ellipticrmk} 
All of the moves described here preserve ellipticity of elements of
$G$. In many cases this can be seen by noting that the move defines an
equivariant map from $X$ to the resulting tree. The case of an expansion
move is discussed in \ref{expansion} below. The only remaining case, the
slide move, follows from the case of an expansion.  As for hyperbolicity,
we note that folds and multi-folds may change hyperbolic elements into
elliptic elements.
\end{remark}

\begin{subsec}[Collapse moves] \label{collapse} 
Let $e\in E(X)$ be an edge such that $G_e = G_{\partial_0 e}$,
and whose endpoints are in different $G$--orbits. To perform a
\emph{collapse move} one simply collapses $e$ and all of its translates
$\gamma e$ (for $\gamma \in G$) to vertices. That is, one deletes
$e$ and identifies its endpoints to a single vertex, and does the same
with translates. The image vertex of $\partial_0 e$ and $\partial_1 e$
will then have stabilizer $G_{\partial_1 e}$. 

The resulting graph $Y$ clearly admits a $G$--action. To see that it is a
tree, let $(e_1, \ldots, e_n)$ be an oriented path in $Y$ without
reversals. The corresponding sequence of edges in $X$ forms a disjoint
union of oriented paths, each without reversals. The unique path in $X$
from $\partial_0 e_1$ to $\partial_1 e_n$ alternates between these paths
and paths in $(G e \cup G \overline{e})$. In particular it is not
contained in $(G e \cup G \overline{e})$, and so the endpoints
$\partial_0 e_1$ and $\partial_1 e_n$ map to different vertices of
$Y$. Thus, the original path in $Y$ is not a circuit. 

To perform a collapse move in ${\bf A}$ one selects an edge $e\in E(A)$
which is not a loop, such that
$\alpha_e\co \As_e \rightarrow \As_{\partial_0 e}$ is an isomorphism. Then
one removes $e$ and $\partial_0 e$, leaving $\partial_1 e$. Every edge
$f$ with initial vertex $\partial_0 e$ is given the new initial vertex
$\partial_1 e$, and each inclusion $\alpha_f$ is replaced with
$\alpha_{\overline{e}} \circ \alpha_e^{-1} \! \circ \alpha_f$. 

A collapse move simplifies the underlying graph of ${\bf A}$ without
increasing any of the labels $\As_v$, $\As_e$. Recall from Definition
\ref{reduced} that ${\bf A}$ is reduced if and only if no collapse
moves can be performed. 
\end{subsec}

\begin{subsec}[Expansion moves] \label{expansion} 
An \emph{expansion move} is the reverse of a collapse move. To perform it
one chooses a vertex $v\in V(X)$, a subgroup $H \subseteq G_v$, and
a collection $S\subseteq E_0(v)$ of edges such that $G_f \subseteq
H$ for every $f \in S$. One then adds a new edge $e$ with $\partial_0 e
= v$, and detaches the edges of $S$ from $v$ and re-defines
$\partial_0 f = \partial_1 e$ for every $f\in S$. Finally, one performs
this operation at each translate $\gamma v$ using the subgroup
$H^{\gamma} \subseteq G_{\gamma v}$ and the edges $\gamma S
\subseteq E_0(\gamma v)$. 

After the move, the stabilizer $G_v$ is unchanged and the new edge
$e$ has stabilizer $H \subseteq G_v$. The stabilizers $G_f$,
for $f \in S$, are also unchanged. Thus the set of elliptic elements 
$\Es = \bigcup_{v \in V(X)} \, G_v$ is the same before and after the
move. 

In ${\bf A}$, an expansion can be performed by choosing a vertex $v\in
V(A)$, a subgroup $H \subseteq \As_v$, and a collection $S \subseteq
E_0(v)$ of edges such that $\alpha_f(\As_f) \subseteq H$ for every $f\in
S$. One then adds a new edge $e$ with $\partial_0 e = v$ and sets $\As_e
= \As_{\partial_1 e} = H$, $\alpha_e = \mbox{ inclusion}$, and
$\alpha_{\overline{e}} = \id$. Finally one
detaches $f$ from $v$ for each $f\in S$ and re-defines $\partial_0 f$ to
be $\partial_1 e$. The inclusions $\alpha_f\co \As_f \rightarrow \As_v$
are unchanged except that they are now regarded as maps into $H$. 
\end{subsec}

\begin{definition}\label{elemdef} 
Collapse moves and expansion moves are called \emph{elementary
moves}. A move which factors as a finite composition of elementary moves
is called an \emph{elementary deformation}. 
Thus slides, subdivision, and the reverse of subdivision are all
elementary deformations (see \ref{slides}, \ref{subdiv} below). 
\end{definition}

\begin{remark} \label{elemsubgroups} 
We have just seen that collapse and expansion moves preserve ellipticity
and hyperbolicity of elements of $G$. In fact a slightly stronger
statement holds: these moves do not change the set of elliptic
\emph{subgroups} of $G$. Recall that a subgroup $H \subseteq G$
is elliptic if and only if it is contained in a vertex stabilizer. During
an elementary move there are two stabilizers involved (up to
conjugacy), and the larger of these two is present before and after the
move. Thus, if $H$ is contained in a stabilizer before the move, it is
still contained in a stabilizer afterward. 
\end{remark}

\begin{subsec}[Slide moves] \label{slides} 
Suppose $e,f\in E(X)$ are adjacent edges (ie, $\partial_0 e =
\partial_0 f$) such that $G_f \subseteq G_e$ and
$f \not\in G e \cup G\overline{e}$. To perform a \emph{slide
move} of $f$ over $e$, detach $f$ from $\partial_0 e$ and re-define
$\partial_0 f$ to be $\partial_1 e$; also do the same for all pairs
$\gamma e, \gamma f$. The requirement that $G_f \subseteq G_e$
ensures that the equivariant move is well defined, and that $G$
still acts on the resulting tree. 

To perform a slide move in ${\bf A}$ one selects adjacent edges $e,f\in
E(A)$ such that $\alpha_f(\As_f) \subseteq \alpha_e(\As_e)$ and $f \not=
e,\overline{e}$. Again one detaches $f$ from $\partial_0 e$ and
re-defines $\partial_0 f$ to be $\partial_1 e$. The inclusion map
$\alpha_f$ is then replaced by $\alpha_{\overline{e}} \circ \alpha_e^{-1}
\circ \alpha_f$. 

A slide move is equal to the composition of an expansion and a
collapse. To see this, consider a slide move (in ${\bf A}$) of $f$ over
$e$. We have that $\alpha_f(\As_f) \subseteq \alpha_e(\As_e) \subseteq
\As_{\partial_0 e}$ and $f \not= e, \overline{e}$. Expand at $\partial_0
e$ using the subgroup $\As_e \subseteq \As_{\partial_0 e}$, pulling
across $e$ and $f$ to the new vertex. Then collapse $e$. The newly
created edge takes the place of $e$ (and has label $\As_e$) and $f$ has
been slid across it. 
\end{subsec}

\begin{subsec}[Subdivision] \label{subdiv} 
This move is straightforward: simply insert a vertex in the interior of
an edge $e$ and do the same for each of the translates $\gamma e$. The
new vertex will have stabilizer $G_e$, as will the two ``halves'' of
$e$. 

In ${\bf A}$, insert a vertex in the interior of some $e \in E(A)$,
and give it the label $\As_e$. Also give the two adjacent ``half-edges''
the label $\As_e$ and let the new inclusion maps be the identity. Note
that a subdivision is a special case of an expansion move (just collapse
either of the new half-edges to undo). 
\end{subsec}

\begin{subsec}[Folds] \label{fold} 
To perform a \emph{fold} one chooses edges $e, f\in E(X)$ with
$\partial_0 e = \partial_0 f$ and identifies $e$ and $f$ to a single
edge ($\partial_1 e$ and $\partial_1 f$ are also identified). 
One also identifies $\gamma e$ with
$\gamma f$ for every $\gamma \in G$, so the resulting graph is an
equivariant quotient space of $X$. It is not difficult to show that the
image graph is a tree. 

The effect of a fold on the quotient graph of groups can vary, depending
on how the edges and vertices involved meet the various 
$G$--orbits. The possibilities are discussed in some detail in
\cite{bestvina:accessibility}. To summarize, the fold is of type B if $e$
or $f$ projects to a loop in $G \backslash X$, and type A otherwise;
and it is of type I, II, or III, accordingly as $V(G \backslash X)$
and $E(G \backslash X)$ both decrease, both remain unchanged, or
only $E(G \backslash X)$ decreases. 

A type B fold is equal to the composition of a subdivision, two type A
folds, and the reverse of a subdivision. Thus in many situations it
suffices to consider only type A folds. 
\end{subsec}

\begin{subsec}[Multi-folds] 
This move is similar to a fold except that several edges are folded
together at once. To perform a multi-fold, one chooses a vertex $v \in
V(X)$ and a subset $S \subseteq E_0(v)$, and identifies all the edges of
$S$ to a single edge. One also does the same for all sets of edges
$\gamma S \subseteq E_0(\gamma v)$ for $\gamma \in G$, so the move is
equivariant. Multi-folds are classified into types A and B, and into
types I, II, and III in the same way that folds are.

If $G \backslash X$ is finite then every multi-fold is a finite
composition of type I and type III folds and type II multi-folds. 

In ${\bf A}$, a type II multi-fold corresponds to ``pulling a subgroup
across an edge,'' described briefly in \cite{bestvina:accessibility}. If
the subgroup is finitely generated then the type II multi-fold is a
finite composition of type II folds. 
\end{subsec}

\begin{subsec}[Parabolic folds] \label{parafold} 
This move may also be viewed as infinitely many folds (or multi-folds)
performed in one step. It is similar to a type II multi-fold
except that it is performed at an end $\varepsilon \in \partial X$ rather
than at a vertex. During the move, various rays with end $\varepsilon$
are identified with each other. The exact definition is fairly technical
so we begin with two examples. 
\end{subsec} 

\begin{subexample}[Example] \label{parafoldA} 
Let $G$ be the Baumslag--Solitar group $BS(1,6)$ with presentation
$\langle \, x,t \mid 
txt^{-1} = x^6 \, \rangle$. This group admits a graph of groups
decomposition in which the graph is a single loop, the vertex and edge
groups are both ${\Z}$, and the inclusion maps are the identity and
multiplication by $6$. The Bass--Serre tree $X$ has vertex 
stabilizers equal to the conjugates of the infinite cyclic subgroup
$\langle x \rangle$. Note that as the tree is parabolic, the set of
elliptic elements $\Es = \bigcup_{i \in {\Z}} \, \langle t^{i} x t^{-i}
\rangle$ is a subgroup. 

Now let $Y$ be the linear ${\Z}$--tree where a generator of ${\Z}$
acts by a translation of amplitude $1$. There is a surjective
homomorphism $BS(1,6) \rightarrow {\Z}$ defined by sending $x$ to
$0$ and $t$ to $1$. This homomorphism, with kernel $\Es$, defines an
action of $BS(1,6)$ on $Y$ in which every vertex and edge stabilizer is
$\Es$. Thus $X$ and $Y$ have the same elliptic and hyperbolic elements
(as $BS(1,6)$--trees). 

As a group, $\Es$ is isomorphic to the additive subgroup ${\Z}[1/6]
\subseteq {\Q}$ generated by integral powers of $6$. 
The quotient graph of groups of $Y$ has a loop as its underlying graph,
and its vertex and edge groups are both ${\Z}[1/6]$. The inclusion maps
are the identity and multiplication by $6$ (which is an isomorphism). 

The tree $Y$ is simply the quotient $G$--tree $\Es
\backslash X$. This transition from $X$ to $Y$, which is a special case
of a parabolic fold, can also be achieved with an infinite sequence of
type II folds which ``zip'' the parabolic tree down to a line. 
\end{subexample} 

\begin{subexample}[Example] \label{parafoldB} 
Again consider $BS(1,6)$ and its Bass--Serre tree $X$. Let $H$ be the
subgroup ${\Z}[1/3] \subseteq {\Z}[1/6] = \Es$. Let $e\in E(X)$ be
an edge with $G_e = {\Z} \subseteq {\Z}[1/6]$, so that
$G_e \subseteq H \subseteq \Es$. Form a quotient space $Y$ of $X$ by
identifying the orbit $H e$ to a single edge, and by extending these
identifications equivariantly. Thus the orbits $(H^{\gamma}) \gamma e$
become edges, for each $\gamma \in G$. 

The image of $e$ in $Y$ has stabilizer equal to $H$, and the quotient
graph is again a loop. The quotient graph of groups has vertex and edge
groups equal to ${\Z}[1/3]$ and the inclusion maps are the identity and
multiplication by $6$. The tree $Y$ is the regular tree of valence
three (cf Remark \ref{parafoldrmk}). 
\end{subexample}

\begin{subexample}[Definition] 
Now we discuss the parabolic fold move. Let $e\in E(X)$ be an edge with
$G_e = G_{\partial_0 e}$ such that $\partial_1 e = t \, \partial_0 e$ for
some $t\in G$.  Let $\varepsilon$ be the end represented by the ray $(e,
te, t^2e, \cdots)$. Note that $t$ is hyperbolic (by Lemma
\ref{hypcriteria}(\ref{h1})) and it fixes $\varepsilon$. Let $T$ be the
connected component of the orbit $G \{e, \overline{e}\}$ that contains
$e$, and let $G_T$ be the stabilizer of $T$. Note that $\gamma \in G_T$
if and only if $[e,\gamma e] \subseteq T$. In particular $G_{\varepsilon}
\subseteq G_T$, since if $\gamma \in G_{\varepsilon}$ then $[e,\gamma e]$
is contained in $[e, \varepsilon) \cup [\gamma e, \varepsilon) \subseteq
T$. 

We claim that $G_T$ fixes $\varepsilon$. To see this, suppose that
$\gamma \in G_T$ does not fix $\varepsilon$. Then there is a linear
subtree $(\varepsilon, \gamma \varepsilon) \subseteq T$. Note that $t^i e
\in (\varepsilon, \gamma \varepsilon)$ for sufficiently large $i$ and
$\partial_0 t^i e$ separates $t^i e$ from $\gamma
\varepsilon$. Similarly, for large $i$, $\gamma t^i e \in (\varepsilon,
\gamma \varepsilon)$ and $\partial_0 \gamma t^i e$ separates $\gamma t^i
e$ from $\varepsilon$. These observations imply that there is a vertex $v
\in (\varepsilon, \gamma \varepsilon)$ such that both edges of $E_0(v)
\cap (\varepsilon, \gamma \varepsilon)$ are in $G e$ (rather than $G
\overline{e}$). Then there is an element $\delta \in G$ taking one of
these edges to the other and fixing their initial endpoint $v$. A
conjugate of $\delta$ fixes $\partial_0 e$ without fixing $e$,
contradicting the assumption that $G_e = G_{\partial_0 e}$. 

Thus $G_T \subseteq G_{\varepsilon}$ and so $T$ is a parabolic
$G_T$--tree. In fact we now have that $G_T = G_{\varepsilon}$. Note that
$\varepsilon$ is the unique fixed end of $T$ that is separated from $e$
by $\partial_1 e$. Therefore $\varepsilon$ is well defined without
reference to $t$. Now let $\Es \subseteq G$ be the set of elliptic
elements. 

Choose a subgroup $H \subseteq (G_{\varepsilon} \cap \Es)$ with
$G_e\subseteq H$. Based on $H$ and $e$, we will define an equivalence
relation on $X$ whose quotient space is a $G$--tree. This tree will be
the ``largest'' quotient $G$--tree with the property that the stabilizer
of the image of $e$ contains $H$. The element $t$ is used in the
definition, but the move is independent of the choice of $t$. 

Let $\Sigma$ be the semigroup $\{t^{-i} \mid i \geq 0\}$. 
We define a relation on $E(X)$ as follows: set $e \approx h^{\sigma} e$
for every $h\in H$ and $\sigma\in \Sigma$, and extend
equivariantly: $\gamma e \approx \gamma h^{\sigma} e$ for
every $\gamma \in G$, $h \in H$, and $\sigma \in \Sigma$. This
relation is reflexive and symmetric, but we must pass to its transitive
closure to obtain an equivalence relation. It will be helpful to express
this transitive closure explicitly. Note that $\gamma e \approx \gamma
h_1^{\sigma_1} e \approx \gamma h_1^{\sigma_1} h_2^{\sigma_2} e \approx
\cdots \approx \gamma h_1^{\sigma_1} \cdots h_j^{\sigma_j} e$, so set 
\begin{equation} \label{foldreln} 
\gamma e \sim \gamma h_1^{\sigma_1} \cdots
h_j^{\sigma_j} e 
\end{equation}
where $\gamma \in G$, $h_i \in H$, and $\sigma_i \in \Sigma$ for
each $i$. Let $\langle H^{\Sigma} \rangle$ denote the subgroup of
$G$ generated by the elements $h^{\sigma}$ ($h \in H$,
$\sigma \in \Sigma$). The relation $\sim$ is symmetric and
reflexive, and clearly $e_0 \approx e_1$ implies $e_0 \sim e_1$.
For transitivity, suppose that $e_0 \sim e_1$ and $e_1 \sim e_2$. We can
write $e_0 = \gamma e$ and $e_1 = \gamma h_1^{\sigma_1} \cdots
h_j^{\sigma_j} e$, and also $e_1 = \delta e$ and $e_2 = \delta
k_1^{\tau_1} \cdots k_l^{\tau_l} e$. Then $e = \delta^{-1} \gamma
h_1^{\sigma_1} \cdots h_j^{\sigma_j} e$. This shows that $\delta^{-1}
\gamma \in G_e \langle H^{\Sigma} \rangle = \langle H^{\Sigma}
\rangle$. Here we use the assumption that $G_e \subseteq H$. We now
have 
\begin{equation*}
\begin{split}
e_2 \ &= \ \delta k_1^{\tau_1} \cdots
k_l^{\tau_l} \delta^{-1} \gamma h_1^{\sigma_1} \cdots h_j^{\sigma_j} e \\
&= \ \gamma (\gamma^{-1} \delta) k_1^{\tau_1} \cdots k_l^{\tau_l}
(\delta^{-1} \gamma) h_1^{\sigma_1} \cdots h_j^{\sigma_j} e,
\end{split}
\end{equation*}
showing that $e_0 \sim e_2$, as $(\gamma^{-1} \delta), (\delta^{-1}
\gamma) \in \langle H^{\Sigma} \rangle$. Thus, $\sim$ is the
transitive closure of $\approx$. 

Next we verify that $\sim$ is independent of the choice of $t$. Suppose
that $\partial_1 e = s \partial_0 e$ for some $s\in G$. Then $s^{-1}
t \partial_0 e = \partial_0 e$ and so $s^{-1} t \in G_{\partial_0 e}
= G_e \subseteq H$. Let $k = s^{-1} t$. Notice that for any $h \in
H$ and $i \geq 0$, 
\begin{equation*}
\begin{split}
h^{t^{-i}} \ &= \ (k^{-1} s^{-1}) \cdots (k^{-1}s^{-1}) h (sk) \cdots
(sk) \\
&= \ k^{-1} (k^{-1})^{s^{-1}} (k^{-1})^{s^{-2}} \cdots
(k^{-1})^{s^{-(i-1)}} h^{s^{-i}} k^{s^{-(i-1)}} \cdots k^{s^{-1}} k. 
\end{split}
\end{equation*}
Thus $e_0 \approx_t e_1$ implies $e_0 \sim_s e_1$, where the subscripts
indicate whether $t$ or $s$ has been used. Similarly $e_0 \approx_s e_1$
implies $e_0 \sim_t e_1$, and therefore $\sim_s$ is the same relation as
$\sim_t$. Furthermore, the subgroup $\langle H^{\Sigma} \rangle$ is
independent of the choice of $t$. 

Now let $Y$ be the quotient $G$--graph $X/\!\sim$.  The quotient map $X
\rightarrow Y$ is called a \emph{parabolic fold}. We claim that $Y$ is a
tree.  First consider the quotient $T / \!\sim$. We begin by observing
that every ray of the form $[x, \varepsilon) \subseteq T$ maps
injectively to $T/\!\sim$. Indeed, if two vertices $u, v \in
[x,\varepsilon)$ are related by $\sim$ then there is an element
$\gamma\in (G_{\varepsilon} \cap \Es)$ taking $u$ to $v$. Then $\gamma
([u,\varepsilon)) = [v,\varepsilon)$, and one of these rays contains the
other. Applying Lemma \ref{hypcriteria}(\ref{h3}) to an edge of $[u,
\varepsilon)$ one finds that $\gamma$ is hyperbolic, unless $u=v$. The
latter must occur since $\gamma \in \Es$. 

Now suppose that $e_0 \sim f_0$ where $e_0, f_0 \in E(T)$. To show that
$T/\!\sim$ is a tree it suffices to show that $[e_0,\varepsilon)$ and
$[f_0,\varepsilon)$ have the same image, so that the image of their union
is a ray in $T / \!\sim$. Write $[e_0, \varepsilon) = (e_0, e_1, e_2,
\ldots)$ and $[f_0, \varepsilon) = (f_0, f_1, f_2, \ldots)$ (these are
oriented paths without reversals). We need to show that $e_i \sim f_i$
for each $i$, given that $e_0 \sim f_0$. So suppose that $e_0 = \gamma e$
and $f_0 = \gamma h_1^{\sigma_1} \cdots h_j^{\sigma_j} e$.  Note that
$[e, \varepsilon) = (e, t e, t^2 e, \ldots)$, and so $e_i = \gamma t^i e$
for each $i$.  Then $f_0 = \gamma h_1^{\sigma_1} \cdots h_j^{\sigma_j}
\gamma^{-1} e_0$, and also $\gamma h_1^{\sigma_1} \cdots h_j^{\sigma_j}
\gamma^{-1}$ fixes $\varepsilon$. Thus $\gamma h_1^{\sigma_1} \cdots
h_j^{\sigma_j} \gamma^{-1}([e_0,\varepsilon)) = [f_0, \varepsilon)$ and
so $f_i = \gamma h_1^{\sigma_1} \cdots h_j^{\sigma_j} \gamma^{-1} e_i$
for each $i$. We now have 
\begin{align}
f_i \ &= \ \gamma h_1^{\sigma_1} \cdots h_j^{\sigma_j} \gamma^{-1} e_i
\notag \\ 
&= \ (\gamma t^i) t^{-i} h_1^{\sigma_1} \cdots h_j^{\sigma_j}
t^i e \notag \\
&= \ (\gamma t^i) h_1^{(t^{-i}\sigma_1)} \cdots
h_j^{(t^{-i}\sigma_j)} e. \label{line}
\end{align}
Note that $t^{-i}\sigma_l \in \Sigma$ for each $l \leq j$ because $i
\geq 0$. Thus, writing $f_i$ as above and $e_i$ as $(\gamma
t^i) e$, we see that $e_i \sim f_i$. Therefore $T/\!\sim$ is a tree. 
It follows that $Y$ is also a tree because all identifications occur
in $T$ and its translates. Line \eqref{line} demonstrates
the need for involving the semigroup $\Sigma$ in the definition, since
the exponents $t^{-i}$ are unavoidable if the quotient space is to
be a tree. 

Next we consider stabilizers in $Y$. These are simply the stabilizers of
the equivalence classes in $X$. We have 
\begin{align}
[e] \ &= \ \{\gamma h_1^{\sigma_1} \cdots h_j^{\sigma_j} e \mid \gamma
\in G_e, \ h_i \in H, \ \sigma_i\in \Sigma, \ j\geq 0 \} \notag \\ 
&= \ \{h_1^{\sigma_1} \cdots h_j^{\sigma_j} e \mid h_i \in H, \
\sigma_i\in \Sigma, \ j\geq 0 \} \tag*{ as $G_e \subseteq H$ } 
\end{align}
and hence $G_{[e]} = \langle H^{\Sigma} \rangle$. Similarly
$G_{[\partial_0 e]} = \langle H^{\Sigma} \rangle$, and the stabilizers of
translates of $[e]$ and $[\partial_0 e]$ are obtained by conjugation. If
$x \in (X - (G e \cup G \overline{e} \cup G \partial_0 e))$ then $[x] =
\{x\}$ and $G_{[x]} = G_x$. Notice that the union of all stabilizers is
unchanged, as $\langle H^{\Sigma} \rangle \subseteq (G_{\varepsilon} \cap
\Es)$. Therefore parabolic folds preserve ellipticity and hyperbolicity
of elements of $G$.

One last observation is that the quotient graph does not change during a
parabolic fold; each equivalence class $[x]$ is contained in the orbit
$G x$. 

\end{subexample} 
To summarize: 
\begin{proposition} \label{parafoldsummary} 
Let $e$ be an edge of a $G$--tree such that $G_e = G_{\partial_0 e}$ and
$\partial_1 e = t \partial_0 e$ for some $t\in G$. 
Suppose a parabolic fold is performed at $e$ using the subgroup $H$,
according to relation \eqref{foldreln}. Then the quotient space is a
$G$--tree and the images of $e$ and $\partial_0 e$ have stabilizers
$\langle H^{\Sigma} \rangle$. Outside of $G e \cup G \overline{e} \cup G
\partial_0 e$ the tree and its stabilizer data are unchanged. The move
preserves ellipticity and hyperbolicity of elements of $G$, and induces
an isomorphism of quotient graphs.  \endproof
\end{proposition}

\begin{remark} \label{parafoldrmk} 
In Examples \ref{parafoldA} and \ref{parafoldB} the tree $Y$ cannot be
obtained from $X$ by an elementary deformation. One way to see this is to
note that the trees have different \emph{modular homomorphisms}, whereas
elementary moves preserve this invariant. The modular homomorphism 
$q\co G \rightarrow {\Q}^{\times}$ of a locally finite $G$--tree is
defined by 
\[q(\gamma) \ = \ [V:V\cap V^{\gamma}] \ / \ [V^{\gamma}: V \cap
V^{\gamma}] \, ,\] 
where $V$ is any subgroup of $G$ commensurable with a vertex
stabilizer. In this definition we are using the fact that in locally
finite $G$--trees, vertex stabilizers are commensurable with all of
their conjugates. One can easily check that $q$ is independent of the
choice of $V$. Equivalent definitions and properties of the modular 
homomorphism are given in \cite{bass:treelat}. 

Recall that during an expansion or collapse move there is a vertex
stabilizer that remains unchanged. Taking $V$ to be this stabilizer, one
obtains invariance of $q$ under elementary moves. 

In the two examples of \ref{parafold}, the image of the modular
homomorphism of $X$ is generated by the element $6$. In Example
\ref{parafoldA}, the tree $Y$ has trivial modular homomorphism. In
Example \ref{parafoldB} the image of the modular homomorphism of $Y$ is
generated by the element $\mbox{$[\,{\Z}[1/3] : 6 \, {\Z}[1/3]\,]$}
= 2$. Evidently the modular homomorphism is not invariant under parabolic
folding, or under infinite compositions of elementary moves.
\end{remark} 

\begin{proposition}\label{bddparafold} 
Let $\rho\co X \rightarrow Y$ be a parabolic fold, performed at $e \in
E(X)$. Then $\rho$ is a finite composition of multi-folds if and only if
$G_{\rho(e)}$ fixes a vertex of $X$. When this occurs, the
multi-folds are of type II. 
\end{proposition}

Proposition \ref{bddparafold} will be proved in section
\ref{equivsection}, using the Multi-fold Lemma (\ref{foldlemma}). 

\begin{proposition}\label{multifoldelem} 
If a fold or a type II multi-fold preserves hyperbolicity of elements of
$G$ then it is a finite composition of elementary moves. 
\end{proposition}

\begin{proof}
First we show that the fold must be of type I or type II. Suppose that
edges $e$ and $f$ generate the fold (so $\partial_0 e = \partial_0 f$) and
that it is of type III. This means that $f \not\in G e$ but
$\partial_1 f \in G \partial_1 e$. Then $\partial_1 f = \gamma
\partial_1 e$ for some $\gamma \in G$. If $f \not\in G \overline{e}$ then
$\gamma$ is hyperbolic by Lemma \ref{hypcriteria}(\ref{h2}). Otherwise,
if $f = \delta \overline{e}$ for some $\delta \in G$, then $\delta
\partial_1 e = \partial_0 f$ and $\delta \partial_0 f = \partial_1
f$. Again by Lemma \ref{hypcriteria}(\ref{h2}), $\delta$ is
hyperbolic. We now re-define $\gamma$ to be $\delta^2$, which is
hyperbolic and takes $\partial_1 e$ to $\partial_1 f$. In both cases,
after the fold, the image vertex of $\partial_1 e$ is fixed by
$\gamma$. Hence $\gamma$ becomes elliptic, a contradiction. 

Thus we consider type I folds and type II multi-folds. We can assume that
the folds are of type A, as remarked in \ref{fold}. Consider first a type
I fold. In this situation, $f \not\in (G e \cup G \overline{e})$ and
$\partial_1 f \not\in (G \partial_1 e \cup G \partial_0 e)$.  Since the
stabilizer of the image of $\partial_1 e$ is $\langle G_{\partial_1 e},
G_{\partial_1 f}\rangle$, this subgroup contains no hyperbolic
elements. Thus it fixes an end or vertex $w \in X \cup \partial X$ by
Proposition \ref{titslemma}.  Consider the paths $[w,\partial_1 e]$ and
$[w, \partial_1 f]$. If $e\in [w, \partial_1 e]$ then $G_{\partial_1 e} =
G_e$ as $G_{\partial_1 e}$ fixes $w$. If in addition $f \not\in
[w,\partial_1 f]$ then $f \in [w,\partial_1 e]$, and so $G_{\partial_1 e}
\subseteq G_f$. The same result as the fold can now be achieved by
sliding $e$ over $f$ and then collapsing $e$.

Similarly, if $f \in [w,\partial_1 f]$ and $e \not\in [w, \partial_1 e]$
then the fold is equivalent to sliding $f$ over $e$ and then collapsing
$f$. 

The last possibility is that $e \in [w,\partial_1 e]$ and $f\in [w,
\partial_1 f]$ (the three cases correspond to whether the path from $w$
to $[e,f]$ joins at $\partial_1 f$, at $\partial_1 e$, or at $\partial_0
e = \partial_0 f$). Then $G_e = G_{\partial_1 e}$ and $G_f =
G_{\partial_1 f}$ as before.  In this case we collapse both $e$ and $f$
end expand a new edge with stabilizer $\langle G_e, G_f \rangle$,
bringing across all edges which were previously incident to $\partial_1
e$ or $\partial_1 f$. The result is the same as the fold. 

Finally, consider a type II multi-fold, performed at the set of edges $S
\subseteq E_0(v)$. In this situation $S$ is contained in a
single edge orbit $G e$, with $e\in S$. Consider $G_e$ and
$G_{\partial_1 e}$. Assuming that $\abs{S} > 1$ (otherwise the
multi-fold is trivial) there is an edge $\gamma e \in S$ with $\gamma
\in (G_v - G_e)$. The image vertex of $\partial_1 e$ and
$\partial_1 \gamma e$ has stabilizer containing $G_{\partial_1 e}$
and $G_{\partial_1 \gamma e}$, and therefore $\langle
G_{\partial_1 e}, G_{\partial_1 \gamma e} \rangle$ contains
only elliptic elements. From the Hyperbolic Segment Condition
(\ref{hypsegment}) it follows that $G_{\partial_1 e} = G_e$
(and $G_{\partial_1 \gamma e} = G_{\gamma e}$). The multi-fold
can be replaced by the following: collapse $e$ (the rest of $S$ collapses
with it, by equivariance) and then expand a new edge with the appropriate
stabilizer. 
\end{proof}

\section{Deformation of $G$--trees}\label{equivsection} 

In this section we prove several of the implications comprising Theorem
\ref{maindefthm}. 

\begin{definition}
A group $G$ has property (FA) if every $G$--tree is elliptic. This notion
is due to Serre. It is \emph{unsplittable} if, in every $G$--tree, every
element of $G$ is elliptic. Being unsplittable is equivalent to the
property that $G$ admits no nontrivial graph of groups decomposition. 

These properties often agree, but they do not agree in general. We say
that $G$ has property (E) if there is a global fixed point in every
$G$--tree in which the elements of $G$ are all elliptic. Thus, a group
has property (FA) if and only if it is unsplittable and has property
(E). Note that all finitely generated groups have property (E), by
Proposition \ref{titslemma}. 
\end{definition}

\begin{theorem} \label{mainthm} 
Let $G$ be a group, and let $X$ and $Y$ be cocompact $G$--trees. Then $X$
and $Y$ are related by an elementary deformation if and only if they have
the same elliptic subgroups.
\end{theorem}

\begin{corollary}\label{mainthmcor} 
Let $G$ be a group. Let $X$ and $Y$ be cocompact $G$--trees whose vertex
stabilizers have property \textup{(E)}. The following conditions are
equivalent. 
\begin{enumerate}
\item \label{t1} $X$ and $Y$ are related by an elementary deformation. 
\item \label{t2} $X$ and $Y$ define the same partition of $G$ into
elliptic and hyperbolic elements. 
\item \label{t3} The length functions $\ell_X$ and $\ell_Y$ vanish on the
same elements of $G$. 
\end{enumerate}
\end{corollary}

\begin{proof}
The implications (\ref{t2})$\Leftrightarrow$(\ref{t3}) are trivial, and
(\ref{t1})$\Rightarrow$(\ref{t2}) follows from Theorem \ref{mainthm}. For
(\ref{t2})$\Rightarrow$(\ref{t1}), consider a stabilizer $G_x$ of some
vertex $x\in V(X)$. As a group acting on $Y$ it consists of elliptic
elements, and by property (E) it has a fixed point. Thus every $X$--elliptic
subgroup is $Y$--elliptic, and conversely by symmetry. Now Theorem
\ref{mainthm} yields conclusion (\ref{t1}). 
\end{proof}

\begin{remark}
In Corollary \ref{mainthmcor}, the assumption of (E) vertex stabilizers
is essential. The implication
(\ref{t2}),(\ref{t3})$\Rightarrow$(\ref{t1}) can fail to hold for trees
related by a parabolic fold, as shown in Remark \ref{parafoldrmk}. 
\end{remark}

The following result is the conjecture of Herrlich, slightly modified,
from \cite{herrlich}. 

\begin{corollary} \label{herrlichcor} 
Let ${\bf A}$ and ${\bf B}$ be graphs of groups having finite underlying
graphs, whose vertex groups have property \textup{(FA)}. If ${\bf A}$ and
${\bf B}$ have isomorphic fundamental groups then there exist a graph of
groups ${\bf B'}$, an isomorphism of graphs of groups $\Phi\co {\bf B'}
\to {\bf B}$, and a finite sequence of elementary moves taking ${\bf A}$
to ${\bf B'}$. 
\end{corollary}

The isomorphism of graphs of groups is meant in the sense of \cite[\S
2]{bass:covering}. This notion of isomorphism is more general than the
naive notion (consisting of group isomorphisms satisfying the appropriate
commutative diagrams). For example, without changing the isomorphism
type, one can replace an inclusion map $\alpha_e \co \As_e \to
\As_{\partial_0 e}$ by $\ad(g) \circ \alpha_e \co \As_e \to
\As_{\partial_0 e}$, for any $g \in \As_{\partial_0 e}$. Here, $\ad(g)\co
\As_{\partial_0 e} \to \As_{\partial_0 e}$ is defined by $\ad(g)(s) = g s
g^{-1}$.

\begin{proof}
Write ${\bf A} = (A,\As,\alpha)$ and ${\bf B} = (B, \Bs, \beta)$. Choose
basepoints $a_0 \in A$ and $b_0 \in B$, and let $\psi \co \pi_1({\bf
A},a_0) \to \pi_1({\bf B},b_0)$ be an isomorphism. Setting $G =
\pi_1({\bf A},a_0)$, the trees $X = (\widetilde{{\bf A},a_0})$ and $Y =
(\widetilde{{\bf B},b_0})$ are $G$--trees (via $\psi$ in the case of
$Y$). The vertex stabilizers of both trees have property (FA), so each
stabilizer has fixed points in both trees. Hence $X$ and $Y$ have the
same elliptic subgroups. Theorem \ref{mainthm} now provides a sequence of
elementary moves from $X$ to $Y$. There is a corresponding sequence of
elementary moves of graphs of groups taking ${\bf A}$ to a quotient graph
of groups ${\bf B'}$ of $Y$. Now, since ${\bf B'}$ and ${\bf B}$ have
$G$--isomorphic covering trees, there is an isomorphism $\Phi \co {\bf
B'} \to {\bf B}$ by \cite[4.2--4.5]{bass:covering}.
\end{proof}

Next we turn to the proof of the Theorem \ref{mainthm}, which 
occupies the rest of this section. 

\begin{definition}
A map between trees (or graphs) is a \emph{morphism} if it sends vertices
to vertices and edges to edges (and respects the maps $\partial_0$,
$\partial_1$, $e \mapsto \overline{e}$). Geometrically, it is a simplicial
map which does not send any edge into a vertex. 
\end{definition}

\begin{proposition}\label{bfprop} 
Let $G$ be a group and suppose that $\phi\co X \rightarrow Y$ is an
equivariant morphism of $G$--trees. Assume further that $\phi$ is
surjective and $G \backslash X$ is finite. Then there exist a
$G$--tree $Z$ and equivariant simplicial maps
$X\stackrel{\rho}{\rightarrow} Z \stackrel{\psi}{\rightarrow} Y$ with
$\phi = \psi \rho$, such that $\rho$ is a finite composition of folds and
$\psi$ induces an isomorphism of quotient graphs. 
\endproof 
\end{proposition}

This result is proved in 
\cite[Proposition, \S 2]{bestvina:accessibility}. The hypotheses are
slightly different but their proof still works. The quotient graph
$G \backslash Y$ is finite by equivariance and surjectivity, rather
than by finite generation and minimality. 

\begin{lemma}[Multi-fold Lemma] \label{foldlemma}
Let $\phi\co X\rightarrow Y$ be an equivariant morphism of
$G$--trees. If $T\subseteq X$ is a simplicial subtree of finite
diameter such that $\phi(T)$ is finite, then $\phi$ factors as
$X\stackrel{\rho}{\rightarrow} X' \stackrel{\phi'}{\rightarrow} Y$ where
$\rho$ is a finite composition of multi-folds and $\phi' \vert_{\rho(T)}$
is an embedding. 
\end{lemma}

\begin{proof}
To simplify the argument we 
use \emph{non-equivariant} folds. We will prove that a 
factorization $\phi = \phi' \rho$ exists, where $\rho$ is a finite
composition of non-equivariant multi-folds and $\phi' \vert_{\rho(T)}$ is
an embedding. If one then replaces the multi-folds of $\rho$ by their 
equivariant counterparts, the factorization is still valid. 

We begin by defining an equivalence relation $\sim_{\phi}$ on the set
$E(T) \times E(T)$. Choose a basepoint $v\in V(T)$ and let $[e,e']_v$
denote the subtree spanned by $e$, $e'$, and $v$ (for $e,e' \in E(T)$). 
Set $(e,e') \sim_{\phi} (f,f')$ if there exists a graph isomorphism
$[e,e']_v \rightarrow [f,f']_v$ fixing $v$ and sending $e$ to $f$ and
$e'$ to $f'$, such that the diagram 
\begin{equation}\label{diagram} 
\begin{split} 
\xymatrix@R-=0pt{ 
[e,e']_v \ar[dd]  \ar[drr]^{\phi} && \\
&& {\phi(T)} \\
[f,f']_v \ar[urr]_{\phi} && }
\end{split}
\end{equation}
commutes. We will call the equivalence classes
\emph{$(T,v,\phi)$--types}. Since the tripods $[e,e']_v$ have bounded 
diameter and $\phi(T)$ is finite, there are only finitely many
$(T,v,\phi)$--types. 

A subtree $S \subseteq T$ \emph{supports} a $(T,v,\phi)$--type if the
$(T,v,\phi)$--type contains a pair $(e,e')$ such that $[e,e']_v \subseteq
S$. By taking a union of such tripods one can find a finite subtree $S$
which supports all $(T,v,\phi)$--types. For such an $S$ we claim that if
$\phi \vert_S$ is an embedding, then $\phi \vert_T$ is an embedding. For
otherwise, if $e \not= e'$ and $\phi(e) = \phi(e')$ for some pair of
edges $e, e' \in E(T)$, then the $(T,v,\phi)$--type containing $(e,e')$
is not supported by $S$ (as $\phi \vert_{[e,e']_v}$ is not injective). 

A \emph{full} multi-fold along $e$ is the multi-fold defined by
identifying all of the edges $E_0(\partial_0 e) \cap \phi^{-1} \phi(e)$
to one edge. It is the ``maximal'' multi-fold along $e$ through which
$\phi$ factors. 

Now let $S\subseteq T$ be a finite tree which supports all
$(T,v,\phi)$--types, and suppose that $\phi \vert_S$ is not an embedding.
Then there are two vertices $u,w\in V(S)$ with the same image in $Y$. 
Since $Y$ contains no circuits, the image of the path from $u$ to $w$ must
contain a reversal. Thus there exists a pair of edges $e_0,e_0' \in E(S)$
such that $\partial_0 e_0 = \partial_0 e_0'$ and $\phi(e_0) =
\phi(e_0')$. Among all such pairs choose $e_0$ and $e_0'$ to minimize the
distance $d$ from $\partial_0 e_0$ to $v$. Factor $\phi$ as $\phi = \phi'
\rho$ where $\rho$ is the full multi-fold along $e_0$. We claim that
$\rho$ preserves the equivalence relation, ie, that $(e,e') \sim_{\phi}
(f,f')$ implies $(\rho(e),\rho(e')) \sim_{\phi'} (\rho(f),\rho(f'))$. 

To show this, let $T_d \subseteq T$ be the subtree spanned by the
vertices having 
distance at most $d$ from $v$, and define $S_d \subseteq S$ similarly
(note, $S_d \subseteq T_d$). By the choice of $e_0$, $\phi \vert_{S_d}$
is an embedding. In fact, since $S$ supports all $(T,v,\phi)$--types,
$\phi \vert_{T_d}$ is an embedding, and then we must have $S_d =
T_d$. This last statement holds because if there were an edge $e\in
E(T_d) - E(S_d)$ then by injectivity of $\phi \vert_{T_d}$, the
$(T,v,\phi)$--type containing $(e,e)$ would not be supported by $S$. 

Now suppose that $(e,e') \sim_{\phi} (f,f')$ and consider the effect of
$\rho$ on the tripods $[e,e']_v$ and $[f,f']_v$. Let $\alpha\co [e,e']_v
\rightarrow [f,f']_v$ be the graph isomorphism shown in diagram
\eqref{diagram}. Note that since $\phi \vert_{T_d}$ is injective, 
the restriction $\alpha \vert_{[e,e']_v \cap T_d}$ is the
identity (and $[e,e']_v \cap T_d \ = \ [f,f']_v \cap T_d$). If
$\partial_0 e_0 \not\in [e,e']_v$ then $\partial_0 e_0 \not\in
[f,f']_v$ and $\rho$ is injective on both trees (here we use
non-equivariance of $\rho$). Therefore $(\rho(e),\rho(e')) \sim_{\phi'}
(\rho(f),\rho(f'))$ via the isomorphism 
\[(\rho \vert_{[f,f']_v}) \ \alpha
\ (\rho \vert_{[e,e']_v})^{-1} \co [\rho(e), \rho(e')]_{\rho(v)} \rightarrow
[\rho(f), \rho(f')]_{\rho(v)}\] 
(cf diagram \ref{diagram2} below). 
If $\partial_0 e_0 \in [e,e']_v$ then $\partial_0 e_0 \in [f,f']_v$ and
$\alpha$ maps $[e,e']_v \cap (E_0(\partial_0 e_0) \cap \phi^{-1}
\phi(e_0))$ bijectively to $[f,f']_v \cap (E_0(\partial_0 e_0) \cap
\phi^{-1} \phi(e_0))$. Note that $\rho$ identifies $E_0(\partial_0 e_0)
\cap \phi^{-1} \phi(e_0)$ to a single edge and is bijective on other
geometric edges. Hence two edges of $[e,e']_v$ are identified
by $\rho$ if and only if their images under $\alpha$ are identified by
$\rho$. This implies that there is a graph isomorphism $\alpha'$ making
the following diagram commute:
\begin{equation} \label{diagram2}
\begin{split} 
\xymatrix@R-=0pt{ 
[e,e']_v \ar[dd]^{\alpha} \ar[rr]^{\rho} &&
\rho([e,e']_v) \ar[dd]^{\alpha'} \ar[drr]^{\phi'}
&& \\ 
&&&& {\phi'(\rho(T))} \\
[f,f']_v \ar[rr]^{\rho} &&\rho([f,f']_v)
\ar[urr]_{\phi'} && } 
\end{split}
\end{equation}
Now restrict $\alpha'$ to the subtree $[\rho(e),\rho(e')]_{\rho(v)}
\subseteq \rho([e,e']_v)$ to obtain a graph isomorphism $[\rho(e),
\rho(e')]_{\rho(v)} \rightarrow [\rho(f),\rho(f')]_{\rho(v)}$ which
realizes the equivalence $(\rho(e),\rho(e')) \sim_{\phi'}
(\rho(f),\rho(f'))$. Thus, $\rho$ preserves the equivalence relation. 

It follows that $\rho(S)$ supports all $(\rho(T), \rho(v), \phi')$--types.
Note also that $\rho(S)$ has fewer edges than $S$. By repeating this
factorization process at most $\abs{E(S)}$ times, we obtain a finite
composition of multi-folds $\rho$ and a morphism $\phi'$ satisfying $\phi =
\phi' \rho$, such that $\phi'\vert_{\rho(S)}$ is an embedding and
$\rho(S)$ supports all $(\rho(T),\rho(v),\phi')$--types. Then
$\phi'\vert_{\rho(T)}$ is an embedding. 
\end{proof}

\begin{proposition} \label{preimage} 
Let $\psi\co Z \rightarrow Y$ be an equivariant morphism of
$G$--trees which induces an isomorphism of quotient graphs, and
which preserves hyperbolicity of elements of $G$. 
\begin{enumerate}
\item \label{p1} If $e$ is any edge of $Y$ with $G_{\partial_0 e}
\supsetneq G_e \subsetneq G_{\partial_1 e}$ then
$\psi^{-1}(e)$ consists of a single edge. 
\item \label{p2} Assume further that $Y$ is reduced and 
$G \backslash Y$ has finitely many vertices. Then 
$\psi$ is a finite composition of parabolic folds. 
\end{enumerate}
\end{proposition}

Note that if every edge of $Y$ satisfies the condition in (\ref{p1}) then
\ref{preimage}(\ref{p1}) implies that $\psi$ is an embedding. Then
$\psi$ will be an isomorphism by equivariance and the assumption on
quotient graphs. 

\begin{proof}[Proof of \ref{preimage}(\ref{p1})] 
Suppose $Z$ contains edges $f$ and $f'$ with $\psi(f) =
\psi(f') = e$. First we reduce to the case where $\partial_0 f =
\partial_0 f'$. Apply the Multi-fold Lemma (\ref{foldlemma}) to the
subtree $[f,f']$ to obtain a factorization of $\psi$:  
\[ Z  \stackrel{\rho_1}{\longrightarrow}  Z_1 
\stackrel{\rho_2}{\longrightarrow} \ \cdots \
\stackrel{\rho_{k-1}}{\longrightarrow}  Z_{k-1}  
\stackrel{\rho_k}{\longrightarrow}  Z_k  \longrightarrow  Y \] 
where each $\rho_i$ is a multi-fold and $\rho_k \cdots \rho_1 (f) =
\rho_k \cdots \rho_1 (f')$. 
Each map in this factorization induces an isomorphism of quotient graphs,
because their composition does. Similarly, each map preserves
hyperbolicity of elements of $G$, because equivariant maps cannot make
elliptic elements hyperbolic. 
Take $i$ smallest so that $\rho_i \cdots
\rho_1 (f) = \rho_i \cdots \rho_1 (f')$. Then 
$\partial_0(\rho_{i-1} \cdots \rho_1 (f)) = \partial_0(\rho_{i-1} \cdots
\rho_1 (f'))$ or $\partial_0(\rho_{i-1} \cdots \rho_1 (\overline{f})) =
\partial_0(\rho_{i-1} \cdots \rho_1 (\overline{f}'))$. We can now proceed
with the argument using this pair of edges and the morphism $Z_{i-1}
\rightarrow Y$. We have already remarked that the hypotheses of the
proposition hold for the this map. 

Thus, without loss of generality, $\partial_0 f = \partial_0 f'$. Since
$\psi(f) = \psi(f')$ and $\psi$ induces an isomorphism on quotient
graphs, $f$ and $f'$ are in the same orbit, so there exists $\gamma\in
G$ with $\gamma f = f'$. Note that $\gamma \in
G_{\partial_0 f} - G_f$, and also $\gamma \in G_e$. Now
choose an element $\sigma \in G_{\partial_1 e} - G_e$. Suppose
that $\sigma$ is in $G_{\partial_1 f} - G_f$ as well. Then
using the Hyperbolic Segment Condition (\ref{hypsegment}) one finds that
$\sigma \gamma$ is hyperbolic in $Z$. However, $\sigma$ and $\gamma$ are
both in $G_{\partial_1 e}$ so their product is elliptic in $Y$,
giving the desired contradiction. Thus our strategy will be to pin down
the behavior of $\sigma$ acting on $Z$. 

Choose an element $\delta \in G_{\partial_0 e} - G_e$. The
product $\sigma \delta$ is hyperbolic in $Y$ by the Hyperbolic Segment
Condition (\ref{hypsegment}) and hence is hyperbolic in $Z$ by
equivariance. Consider its axis 
$Z_{\sigma \delta} \subseteq Z$. The image $\psi(Z_{\sigma \delta})$ is a
$\langle \sigma \delta\rangle $--invariant subtree of $Y$ and therefore
it contains the axis $Y_{\sigma \delta}$. Notice that $Y_{\sigma \delta}$
is the support of the following oriented path: 
\[\cdots \ \ \delta^{-1} \overline{e}, \ \ e, \ \ \sigma \overline{e},
\ \ \sigma \delta e \ \ \cdots .\] 
Since $e\in Y_{\sigma \delta}$, the axis
$Z_{\sigma \delta}$ contains an edge which is mapped by $\psi$ to $e$,
and this edge must be in the same orbit as $f$; 
call it $\tau f$ where $\tau\in G$. Now, using equivariance, one sees
that the set $Z_{\sigma \delta} \cap 
\psi^{-1}(Y_{\sigma \delta})$ contains the following oriented edges,
arranged coherently and in order along $Z_{\sigma \delta}$ (possibly with
gaps in between): 
\begin{equation}\label{axis} 
\cdots \ \ \delta^{-1} \tau \overline{f}, \ \ \tau f, \ \ \sigma
\tau \overline{f}, \ \ \sigma \delta \tau f \ \ \cdots .
\end{equation}
This shows that the edges $\tau f$ and $\sigma \tau f$ are incoherently
oriented, and so 
$\sigma$ fixes a unique vertex $v \in [\tau f, \sigma \tau f]
\subseteq Z_{\sigma \delta}$. Configuration \eqref{axis} also shows that
$[\tau f, \sigma \tau f] = [\partial_0(\tau f), \partial_0(\sigma \tau
f)]$. Hence $[\partial_0(\tau f), v] = [\tau f,v]$, and $\sigma$ fixes no
edge of this segment. 

Having established the behavior of $\sigma$ we return our attention to
$\gamma$, or rather to its conjugate $\tau \gamma \tau^{-1}$. This
element fixes $\partial_0(\tau f)$ and takes $\tau f$ to $\tau f'$
($\not= \tau f$). Now $\sigma$ and $\tau \gamma \tau^{-1}$ satisfy the
Hyperbolic Segment Condition (\ref{hypsegment}) using the segment $[\tau
f, v]$. This implies that $\sigma \tau \gamma \tau^{-1}$ is hyperbolic in
$Z$.  On the other hand, as $\psi(\tau f) = \psi(f) = e$, we must have
$\tau \in G_e$ by equivariance. Then $\sigma$, $\tau$, and $\gamma$ are
in the subgroup $G_{\partial_1 e}$ and so $\sigma \tau \gamma \tau^{-1}$
is elliptic in $Y$, a contradiction.
\end{proof}

\begin{proof}[Proof of \ref{preimage}(\ref{p2})] 
Suppose first that $G \backslash Y$ is a single loop. If the edges
of $Y$ satisfy $G_e \not= G_{\partial_0 e}$ then by
\ref{preimage}(\ref{p1}) the map $\psi$ is an isomorphism. Otherwise
suppose that $G_f = G_{\partial_0 f}$ for some $f\in E(Y)$. Then $Y$
is parabolic and hence the elliptic elements form a subgroup of $G$
(see Remark \ref{parabolicsubgroup}). As 
$G \backslash Z$ is also a single loop, this implies that $Z$ is 
also parabolic (using the Hyperbolic Segment Condition
(\ref{hypsegment})). Let $e\in E(Z)$ be an edge with image $f$ or 
$\overline{f}$, such that $G_e = G_{\partial_0 e}$ and
$\partial_1 e = t \partial_0 e$ for some $t \in G$ (for every edge
$e \in E(Z)$, either $e$ or $\overline{e}$ has these last two
properties, because $Z$ is parabolic). 

We wish to know that $G_{\psi(e)} = G_{\partial_0
\psi(e)}$. If $\psi(e) = f$ then this is already true. So suppose that
$\psi(e) = \overline{f}$. By equivariance, we have that $\partial_1
\overline{f} = t \partial_0 \overline{f}$, ie, $\partial_0 f = t
\partial_1 f$. Hence $G_{t f} \subseteq G_{\partial_1 t f} =
G_{\partial_0 f} = G_f$, which implies that $G_{t^{i+1} f}
\subseteq G_{t^i f}$ for all $i$. Similarly, $G_{t^i e}
\subseteq G_{t^{i+1} e}$ for all $i$, because $G_e \subseteq
G_{\partial_1 e} = G_{t \partial_0 e} = G_{t e}$. These
inclusions fit into a diagram as follows. 
\setcounter{MaxMatrixCols}{13}
\begin{equation}
\begin{matrix}
&&\cdots & \subseteq & G_{t^{-1} e} & \subseteq & G_e &
\subseteq & G_{t e} & \subseteq & \cdots & \subseteq & \bigcup_n
G_{t^n e} \\
&&&&\cap&&\cap &&\cap &&&&\cap \\
\bigcup_n G_{t^n {f}} & \supseteq & \cdots & \supseteq &
G_{t^{-1} {f}} & \supseteq & G_{{f}}&
\supseteq & G_{t {f}} & 
\supseteq & \cdots & \supseteq & \bigcap_n G_{t^n {f}}
\end{matrix}\tag*{}
\end{equation}
The central vertical inclusions follow from equivariance of $\psi$, since
$\psi(t^i e) = t^i \overline{f}$ and $G_{t^i \overline{f}} =
G_{t^i f}$. 
The rightmost vertical inclusion is a consequence of all of the other
inclusions in the diagram. Now observe that $\bigcup_n G_{t^n {f}}$
and $\bigcup_n G_{t^n e}$ are both equal to $\Es$, the set of
elliptic elements. In each case this follows from the facts that $t$ is
hyperbolic (by Lemma \ref{hypcriteria}) and the tree is parabolic. 
Therefore all of the inclusions in the bottom row of the diagram are
equalities. This implies that $Y$ is a linear tree 
with all stabilizers equal to $\Es$, as in Example \ref{parafoldA}. In
particular, $G_{f} = G_{\partial_1 {f}}$ as desired.

Renaming $\overline{f}$ as $f$ if necessary, we have: $\psi(e) = f$,
$G_e = G_{\partial_0 e}$, $G_f = G_{\partial_0 f}$,
$\partial_1 e = t \partial_0 e$, and $\partial_1 f = t \partial_0 f$. Now
consider the subgroup $H = G_f$, and let $\sim$ be the relation on
$X$ defined by \eqref{foldreln} using $H$, $e$, and $t$. We claim that 
\begin{equation} \label{relns} 
e_0 \sim e_1 \ \ \text{ if and only if } \ \ \psi(e_0) =
\psi(e_1).
\end{equation}
This equivalence implies that $\psi$ is a parabolic fold from $Z$ to its
image. Since $Y$ is minimal (because it is reduced), this image is all of
$Y$. 

Note that both relations in \eqref{relns} are equivariant, so in its
proof we can take one of the edges to be $e$. So suppose that $\psi(e') =
\psi(e)$ ( $= f$). Then $e' = \gamma e$ for some $\gamma \in G$, and
$\gamma \in H$ by equivariance. Therefore $e\sim e'$. 

Next suppose that $e \approx e'$, ie, that $e' = h^{t^{-i}} e$ with $i
\geq 0$ and $h \in H$. Then $\psi(e') = h^{t^{-i}} \psi(e) = h^{t^{-i}}
f$. The properties $G_f = G_{\partial_0 f}$ and $\partial_1 f =
t \partial_0 f$ imply that $G_{t^j f} \subseteq G_{t^k f}$
whenever $j \leq k$. Thus $(G_f)^{t^{-i}} = G_{t^{-i} f}
\subseteq G_f$, and so $h^{t^{-i}} \in G_f$. This implies that
$\psi(e') = f = \psi(e)$, and so we have shown that $e_0 \approx e_1$
implies $\psi(e_0) = \psi(e_1)$. Since the latter relation is transitive
and $\sim$ is the transitive closure of $\approx$, it now follows that
\eqref{relns} holds. Therefore \ref{preimage}(\ref{p2}) is proved in the
case where the quotient graph is a single loop. 

Next consider the general case. 
An edge $e\in E(Z)$ is called \emph{$\psi$--isolated} if no two edges of
$\psi^{-1} \psi (e)$ have a common endpoint. 
As $Y$ is reduced, \ref{preimage}(\ref{p1}) implies that every edge that
is not $\psi$--isolated projects to a loop in $G \backslash Z$. We are
interested in $\psi$--isolated edges because of the following
observation: if $\psi \co Z \to Y$ is any morphism of trees such that every
edge of $Z$ is $\psi$--isolated, then $\psi$ is an embedding. This is true
because otherwise, $\psi$ would factor through a nontrivial fold, and
then the edge where the fold is performed would not be $\psi$--isolated. 

Suppose that $e\in
E(Z)$ is an edge which is not $\psi$--isolated. Let
$f = \psi(e)$, and let $e$ be oriented so that $G_e =
G_{\partial_0 e}$. 
Let $T \subseteq Z$ be the connected component of $G e$ containing
$e$. Then $T$ is a parabolic $G_T$--tree with quotient graph a
single loop ($G_T$ is the stabilizer of $T$). Let $S\subseteq Y$ be
the connected component of $G f$ containing $f$. Then $S$ is also a
$G_T$--tree, and $\psi \vert_T \co T \to S$ is a
$G_T$--equivariant map inducing an isomorphism of quotient graphs. By
the previous analysis, $\psi \vert_T$ is a parabolic fold of
$G_T$--trees. In fact it is the parabolic fold performed at $e$ using
the subgroup $H = (G_T)_f$. 

Let $\sim_T$ denote the relation on $T$ given by this parabolic fold, and
let $\sim$ be its equivariant extension to $Z$. Then $\sim$ is simply the
parabolic fold in $Z$ performed at $e$ using $H$. We have that $e_0
\sim_T e_1$ if and only if $\psi(e_0) = \psi(e_1)$, for edges $e_0, e_1
\in E(T)$. Since $\psi$ is equivariant, we also have that $e_0 \sim e_1$
implies $\psi(e_0) = \psi(e_1)$ for all $e_0, e_1 \in E(Z)$. It follows 
that $\psi$ factors as $Z\stackrel{\rho}{\rightarrow} Z'
\stackrel{\psi'}{\rightarrow} Y$, where $\rho$ is the parabolic fold
defined by $\sim$. Note that $\psi' \vert_{\rho(T)} \co \rho(T) \to S$ is
an isomorphism. In particular the parabolic fold $\rho$ is nontrivial,
since $\psi \vert_T \co T \to S$ is not an isomorphism (as $e$ 
is not $\psi$--isolated). 

Recall that in a parabolic fold at $e$, all identifications of vertices
occur inside connected components of $G e$. Therefore edges and vertices
in different components remain in different components. Now let $e' =
\rho(e) \in E(Z')$. We do not know whether ${\psi'}^{-1} \psi'(e') = \{ e'
\}$, but no two edges of ${\psi'}^{-1} \psi'(e')$ are in the same
connected component of $G e'$, and therefore $e'$ is $\psi'$--isolated.

Thus, by factoring through a parabolic fold, any edge may be made
$\psi'$--isolated. During this parabolic fold, the other edges are
affected as follows. All translates of $\rho(e)$ and $\rho(\overline{e})$
are $\psi'$--isolated. If $f$ is a $\psi$--isolated edge of $Z$ such that
$\partial_0 f, \partial_1 f \not\in G \partial_0 e$, then $\rho(f)$ is
$\psi'$--isolated, as $\rho$ does not affect $Z$ away from $G e$. If $f$
is a $\psi$--isolated edge of $Z$ with $\partial_0 f, \partial_1 f \in G
\partial_0 e$ and $f \not\in (G e \cup G \overline{e})$, then
$G_{\partial_0 \rho(f)} \supsetneq G_{\rho(f)} \subsetneq G_{\partial_1
\rho(f)}$. This is true because $G_{\partial_0 \rho(e)}$ is strictly
larger than $G_{\partial_0 e}$, as $\rho$ is nontrivial, and
$G_{\rho(f)} = G_f$ (because $f \not\in (G e \cup G
\overline{e})$). Thus, $\rho(f)$ is $\psi'$--isolated by
\ref{preimage}(\ref{p1}). Lastly, if only one of $\partial_0 f$,
$\partial_1 f$ is in the orbit $G \partial_0 e$, then $\rho(f)$ does not
project to a loop in $G \backslash Z'$, and hence is $\psi'$--isolated by
\ref{preimage}(\ref{p1}).

These remarks imply that $\psi$ factors through a composition of finitely
many nontrivial parabolic folds, at most one for each vertex orbit,
after which all edges are $\psi'$--isolated. Here the factorization is
$Z\stackrel{\rho}{\rightarrow} Z' \stackrel{\psi'}{\rightarrow} Y$, where
$\rho$ now represents the sequence of parabolic folds. Then $\psi'$ is an
embedding, and it is an isomorphism because $Y$ is minimal.
\end{proof}

Notice that the parabolic folds in \ref{preimage}(\ref{p2}) do not
interfere with 
each other and may be performed in any order. Thus a similar result holds
in the case where $G \backslash Y$ has infinitely many vertices, if
one is willing to perform infinitely many ``disjoint'' parabolic folds
simultaneously. 

Next we give a proof of Proposition \ref{bddparafold}, using the
Multi-fold Lemma. Recall that this proposition characterizes parabolic
folds which admit finite factorizations into multi-folds. 

\begin{proof}[Proof of Proposition \ref{bddparafold}] We are given a
parabolic fold $\rho \co X \to Y$, performed at $e \in E(X)$. Suppose
that $G_{\rho(e)}$ fixes a vertex $v\in V(X)$. Let $T \subseteq X$
be the subtree spanned by the edges $\rho^{-1} \rho (e)$. We claim that
$T$ has finite diameter. If $e' \in \rho^{-1} \rho(e)$ then $e' = \gamma
e$ for some $\gamma \in G$ (because parabolic folds induce
isomorphisms on quotient graphs). Then $\gamma \rho(e) = \rho(\gamma e) =
\rho(e') = \rho(e)$, so $\gamma \in G_{\rho(e)} \subseteq
G_v$. We now have $\gamma([v,e]) = [v,e']$ and hence
$d(e,e') \leq 2 d(e,v)$ (as $[e,v] \cup [v,e']$ is a path joining $e$ to
$e'$). This shows that $\rho^{-1} \rho (e)$ has diameter at most
$4 d(e,v)$, and the same must hold for $T$. 

Now, applying the Multi-fold Lemma (\ref{foldlemma}) using the subtree
$T$, we obtain a factorization of $\rho$ as $X
\stackrel{\rho'}{\rightarrow} X' \stackrel{i}{\rightarrow} Y$. Here
$\rho'$ is a finite composition of multi-folds and $i\vert_{\rho'(T)}$ is
an embedding. Thus, $i^{-1}(\rho(e))$ is a single edge. The same is true
of translates of $e$ by equivariance, and also of other geometric edges,
since the original parabolic fold $\rho$ is bijective away from the orbit
of $e$. This shows that $i$ is an isomorphism, and hence $\rho$ has been
factored into a finite composition of multi-folds. 

To see that the multi-folds are of type II, note that their composition
induces an isomorphism on quotient graphs, and hence each multi-fold does
as well. 

Finally consider the converse: suppose that the parabolic fold $\rho$ is a
finite composition of multi-folds (which must be of type II). Note that
these multi-folds occur along edges $f\in G \overline{e}$, so
that $G_f = G_{\partial_1 f}$. It follows that these multi-folds
preserve hyperbolicity of elements of $G$, so by 
Proposition \ref{multifoldelem}, $\rho$ is an elementary deformation. 
Now, by Remark \ref{elemsubgroups}, every stabilizer
in $Y$ fixes a vertex of $X$. In particular, $G_{\rho(e)}$ fixes
some $v\in V(X)$. 
\end{proof}

\begin{lemma}\label{hypconj}\mathsurround = 0pt
Let $X$ be a reduced $G$--tree. 
If $x$ and $y$ are distinct vertices
of $X$ with $G_x \subseteq G_y$ then there is a hyperbolic
element $t\in G$ such that $G_x \subseteq (G_x)^t$. 
\end{lemma}

\begin{proof}
We are given that $G_x$ fixes $y$. Thus $G_x$ fixes the path
from $x$ to $y$, including the first edge $e$ on this path, and so
$G_x = G_e$ (where $\partial_0 e = x$). Since $X$ is reduced,
there is an element $t$ of $G$ taking $\partial_0 e$ to $\partial_1
e$. Then $G_x = G_e \subseteq G_{\partial_1 e} =
(G_x)^t$. The element $t$ is hyperbolic by Lemma \ref{hypcriteria},
using either of the first two criteria. 
\end{proof}

\begin{lemma}\label{hypdist}
Let $Y$ be a $G$--tree and suppose $s,t\in G$ are elements with
$t$ hyperbolic. Let $x$ and $y$ be any two vertices of $Y$. Then $d(x,
st^n y) \rightarrow \infty$ as $n\rightarrow \infty$. 
\end{lemma}

\begin{proof2}
Note that $\ell_Y(t) > 0$ and that $t^{-n}$ has translation length $n \cdot
\ell_Y(t)$ for positive $n$. Thus: 
\begin{align} 
d(x,st^n y) \ = \ d(t^{-n}s^{-1} x,y) \ & \geq
\ d(t^{-n}s^{-1} x,s^{-1} x) - d(s^{-1} x,y) \notag \\
&\geq \ n \cdot \ell_Y(t) - d(s^{-1} x,y). \tag*{\qed} 
\end{align}
\end{proof2}

\begin{proposition} \label{morphism} 
Let $X$ and $Y$ be cocompact $G$--trees which define the same partition of
$G$ into elliptic and hyperbolic elements. Suppose further that $X$ is
reduced and every vertex stabilizer of $X$ fixes a vertex of $Y$. Then,
after possibly subdividing $X$, there is an equivariant morphism $\phi
\co X \to Y$. 
\end{proposition}

\begin{proof} 
We will define the morphism $\phi$ in the following manner. First we
define $\phi$ on a finite set of vertices of $X$, one representing each
vertex orbit. Having done this correctly, the map then extends
equivariantly to $V(X)$ in a unique way. Once the map is defined on
$V(X)$, each edge $e$ is sent to the unique path in $Y$ joining $\phi(
\partial_0 e)$ to $\phi( \partial_1 e)$. We need to be certain that these
two vertices are not the same. Then $e$ is subdivided so that preimages
of vertices in $Y$ are vertices in $X$, and we have an equivariant
morphism $\phi \co X' \to Y$ (here $X'$ is $X$, after subdivision). 

Let $S\subseteq X$ be a subtree whose edges map bijectively to
the edges of $G\backslash X$. For each vertex $v$ of $S$, the 
stabilizer $G_v$ may or may not be maximal among the elliptic
subgroups of $G$. Label the vertices of $S$ as $v_1, \ldots, v_k$ so
that the vertices having maximal elliptic stabilizers occur first in the
list. Also arrange that all vertices in the same orbit occur
consecutively. Consider the $i^{\text{th}}$ orbit arising in the list. It
has a first vertex $v_{i_0}$ and various translates $v_{i_0+q}$ 
for $1 \leq q \leq r_i$. Set $\gamma_{i_0} = 1$ and choose elements
$\gamma_{i_0+q} \in G$ such that $\gamma_{i_0+q} v_{i_0} =
v_{i_0+q}$ for each $q$. 

Thus, whenever $v_i$ and $v_j$ are in the same orbit, there is a chosen
element $\gamma_i \gamma_j^{-1}$ taking $v_j$ to $v_i$. Also, for each
$j$ the element $\gamma_j^{-1}$ takes $v_j$ to a chosen representative of
its orbit. If $i\not= j$ then $\gamma_i \gamma_j^{-1}$ is hyperbolic, by
Lemma \ref{hypcriteria}(\ref{h2}). 

We wish to define $\phi$ on $\{ v_1, \ldots, v_k\}$ with the following
properties: 
\begin{enumerate} 
\renewcommand{\theenumi}{\roman{enumi}} 
\renewcommand{\labelenumi}{(\theenumi)} 
\item \label{d1} $\phi$ is injective, 
\item \label{d2} $G_{v_i} \subseteq G_{\phi(v_i)}$ for each
$i$, 
\item \label{d3} whenever $v_i$ and $v_j$ are in the same orbit,
$\gamma_i \gamma_j^{-1} \phi(v_j) = \phi(v_i)$. 
\end{enumerate}
We proceed one orbit at a time. 
Suppose that $\phi$ has already been defined on $\{v_1, \ldots,
v_{d-1}\}$, satisfying all three properties. Consider $v_d$
and its translates $v_{d+1}, \ldots, v_{d+r}$. The stabilizer
$G_{v_d}$ fixes a vertex $w$ of $Y$. If 
\begin{equation} \label{intersection} 
\{\phi(v_1), \ldots, \phi(v_{d-1})\} \cap \{\gamma_{d} w, \ldots,
\gamma_{d+r} w \} \ = \ \emptyset
\end{equation}
then we can define $\phi(v_{d+q}) = \gamma_{d+q} w$ for each $q\leq
r$. Properties (\ref{d1})--(\ref{d3}) are satisfied in this case.
However, \eqref{intersection} might not hold. So suppose that $\phi(v_p)
= \gamma_{d+q} w$ for some $p \leq d-1$ and $0 \leq q \leq r$. 

We consider two cases: either $G_{v_p}$ 
is maximal among elliptic subgroups, or $G_{v_{d+q}}$ is not
maximal elliptic. The third possibility, that $G_{v_{d+q}}$ is
maximal elliptic and $G_{v_p}$ is not, was ruled out when we 
ordered the vertices. In the first case note that $G_{v_p}
\subseteq G_{\gamma_{d+q} w}$ by property (\ref{d2}) for $v_p$,
and by maximality of $G_{v_p}$ these stabilizers are equal. Also,
since $G_{v_d} \subseteq G_{w}$, conjugation by $\gamma_{d+q}$
yields $G_{v_{d+q}} \subseteq G_{\gamma_{d+q} w}$. Therefore 
$G_{v_{d+q}} \subseteq G_{v_p}$. In the second case, as
$G_{v_{d+q}}$ is not maximal elliptic, there is an elliptic subgroup
$H$ properly containing $G_{v_{d+q}}$. Since $H$ is elliptic it is
contained in $G_v$ for some $v$ and we have $G_{v_{d+q}}
\subsetneq H \subseteq G_v$. 

In each of the two cases we have found that $G_{v_{d+q}}$ is
contained in the stabilizer of a vertex of $X$ other than
$v_{d+q}$. By Lemma \ref{hypconj}, there must be a hyperbolic
element $t'\in G$ such that $G_{v_{d+q}} \subseteq
(G_{v_{d+q}})^{t'}$. Now observe that 
\begin{equation*}
\begin{split}
G_{v_{d}} = (G_{v_{d+q}})^{{\gamma_{d+q}}^{-1}} &\subseteq
((G_{v_{d+q}})^{t'})^{{\gamma_{d+q}}^{-1}} \\ &=
(((G_{v_{d}})^{\gamma_{d+q}})^{t'})^{{\gamma_{d+q}}^{-1}} =
(G_{v_{d}})^{{\gamma_{d+q}}^{-1} t' \gamma_{d+q}}. 
\end{split}
\end{equation*}
We can simplify forthcoming notation by setting $t =
{\gamma_{d+q}}^{-1}t'\gamma_{d+q}$, so that we have $G_{v_{d}} \subseteq
(G_{v_{d}})^t$. This inclusion will allow us to replace $w$ by some $t^n
w$ to arrange that \eqref{intersection} holds. Note that $t$ is
hyperbolic. 

For any pair of indices $p \leq d-1$ and $q \leq r$ consider the vertices
$\phi(v_p)$ and $w$, and the element $s = \gamma_{d+q}$. Applying Lemma
\ref{hypdist} with this data we find that 
\begin{equation} \label{distance} 
d(\phi(v_p), \gamma_{d+q} t^n w) \ > \ 0
\end{equation}
for sufficiently large $n$. Fix $n$ so that \eqref{distance} holds for
every $p \leq d-1$ and $q \leq r$. Now define $\phi(v_{d+q}) =
\gamma_{d+q} t^n w$ for each $q\leq r$. 

Trivially, property (\ref{d3}) holds. For property (\ref{d1}) note first
that since each $\gamma_i \gamma_j^{-1}$ is hyperbolic, $\phi$ is
injective on the set $\{ v_d, \ldots, v_{d+r}\}$. 
It is injective on the set $\{v_1, \ldots, v_{d-1}\}$ by assumption,
and inequality \eqref{distance} implies that the $\phi$--images of these
two sets are disjoint. Thus $\phi$ is injective. Property (\ref{d2})
holds as follows. For each $q\leq r$ we have 
\begin{align}
G_{v_{d+q}} \ &= \ (G_{v_{d}})^{\gamma_{d+q}} \notag \\ 
&\subseteq \ (G_{v_{d}})^{\gamma_{d+q} t^n} \tag*{ as
$G_{v_{d}} \subseteq (G_{v_{d}})^t$ } \\
&\subseteq \ (G_{w})^{\gamma_{d+q} t^n} \tag*{ as 
$G_{v_{d}} \text{ fixes } w$ } \\
&= \ G_{\gamma_{d+q} t^n w} \ = \
G_{\phi(v_{d+q})}. \notag 
\end{align}
Thus $\phi$ satisfies (\ref{d1})--(\ref{d3}), and proceeding inductively
on orbits, such a $\phi$ can be defined on all of $\{v_1, \ldots, v_k\}$. 

Properties (\ref{d2}) and (\ref{d3}) imply that $\phi$ has a (unique)
equivariant extension to $V(X)$.  Now consider any edge $e \in E(S)$,
with endpoints $v_{i_0}$, $v_{i_1}$. Property (\ref{d1}) implies that
$\phi(v_{i_0}) \not= \phi(v_{i_1})$ and so $e$ can be mapped to a
nontrivial embedded path in $Y$. By equivariance this holds for every
edge of $X$, as $S$ meets every edge orbit. Therefore the equivariant
morphism $\phi\co X' \to Y$ can be constructed as required. 
\end{proof}

\begin{proof}[Proof of Theorem \ref{mainthm}] 
We have seen in Remark \ref{elemsubgroups} that elementary moves do not
change the set of elliptic subgroups of $G$. Thus one implication is
clear. For the other implication, we are given cocompact $G$--trees $X$
and $Y$ with the same elliptic subgroups. In particular $X$ and $Y$
define the same partition of $G$ into elliptic and hyperbolic
elements. We can assume without loss of generality that $X$ and $Y$ are
reduced (by performing collapse moves on both).

By Proposition \ref{morphism}, there is an equivariant morphism $\phi \co
X' \to Y$, where $X'$ is the result of subdividing $X$. Note that since
subdivision is a special case of an expansion move, it suffices for us to
show that $X'$ and $Y$ are related by a deformation. As $Y$ is reduced,
it is minimal, and hence $\phi$ is surjective. By Proposition
\ref{bfprop}, $\phi$ factors as $\psi \rho$, where $\psi$ induces an
isomorphism on quotient graphs and $\rho$ is a finite composition of
folds. By Proposition \ref{preimage}(\ref{p2}) the map $\psi$ factors as
a finite composition of parabolic folds, as $\psi = \rho_k \cdots
\rho_0$. Now the original map $\phi$ has the form
\[ X' \stackrel{\rho}{\longrightarrow} Z_0
\stackrel{\rho_0}{\longrightarrow} Z_1 
\stackrel{\rho_1}{\longrightarrow} \cdots 
\stackrel{\rho_{k-1}}{\longrightarrow} Z_{k}
\stackrel{\rho_{k}}{\longrightarrow} Y. \] 
Suppose the parabolic fold $\rho_k$ is performed at $e \in
E(Z_k)$. Consider the stabilizer $G_{\rho_k e}$ in $Y$. By our
assumptions on $X$ and $Y$, this subgroup has a fixed point $v\in
V(X')$. By equivariance of the maps, $G_v$
fixes the vertex $\rho_{k-1} \cdots \rho_0 \rho \, v$ of $Z_k$, and hence
$G_{\rho_k e}$ fixes a vertex of $Z_k$. Now we
can apply Proposition \ref{bddparafold} to the map $\rho_k$ and factor it
into a finite composition of type II multi-folds. Repeating this argument
for each $\rho_i$, we conclude that $\phi$ factors into a finite
composition of folds and type II multi-folds. 

Recall that the map $\phi$ preserves hyperbolicity of elements of
$G$. It follows that each individual fold or multi-fold also has
this property, as equivariant maps cannot make elliptic elements
hyperbolic. Therefore, applying Proposition \ref{multifoldelem} to each
fold and multi-fold, $\phi$ is a finite composition of elementary
moves. This concludes the proof of the theorem. 
\end{proof}

\section{Local Rigidity} \label{localsec} 

\begin{definition}
A $G$--tree is \emph{proper} if $G_e \subsetneq G_{\partial_0 e}$ for
every edge $e$. This is equivalent to the property ``$V(X)$ is a strict
$G$--set'' from \cite{bass:rigidity}. Notice that proper $G$--trees are
reduced and minimal. 

A $G$--tree is \emph{slide-free} if it is minimal and, for all edges $e$
and $f$ with $\partial_0 e = \partial_0 f$, $G_e \subseteq G_f$ implies
$f \in Ge$ or $f \in G \overline{e}$.  A quick look at Definition
\ref{slides} will confirm that a minimal $G$--tree is slide-free if and
only if it admits no slide moves. 

A $G$--tree is \emph{strongly slide-free} if it is minimal and, for all
edges $e$ and $f$ with $\partial_0 e = \partial_0 f$, $G_e \subseteq G_f$
implies $f \in Ge$. 
Consider for example the graph of groups with underlying graph a loop,
with edge and vertex groups ${\Z}$, and with inclusion maps 
multiplication by $2$ and $4$. 
The corresponding tree is slide-free because there is
only one geometric edge orbit. However, it is not strongly slide-free. 
\end{definition}

\begin{remark} \label{ssfremark} 
Strongly slide-free trees are proper, as follows. Let $X$ be a strongly
slide-free $G$--tree and let $e\in E(X)$ be any edge. Minimal trees have
no terminal edges, so there is an edge $e' \not= e$ with $\partial_0 e' =
\partial_0 e$. If $G_e = G_{\partial_0 e}$ then $G_{e'} \subseteq G_e$,
and so $e' = \gamma e$ for some $\gamma\in G$. But then $\gamma \in
G_{\partial_0 e} - G_e$, a contradiction. Therefore $G_e \subsetneq
G_{\partial_0 e}$. 
\end{remark}

Throughout this section, $X$ is a $G$--tree which possesses a maximal
stabilizer. Every proper $G$--tree has this property, for example.  The
set of maximal stabilizers is called $M(X)$. If $V \in M(X)$ then $X_V$
is the tree of fixed points of $V$, called the \emph{characteristic
subtree} of $V$. Note that each $X_V$ is nonempty and $X_V \cap X_W =
\emptyset$ for all distinct pairs $V,W \in M(X)$. Let $[V,W]$ denote the
unique smallest segment joining $X_V$ to $X_W$. Note that $[V,W]$ has
positive length, and meets $X_V$ and $X_W$ only in its endpoints. 

If $G$ contains a hyperbolic element then there is a unique minimal
$G$--invariant subtree $X_G$, which is nonempty (see
\cite[7.5]{bass:covering}). This occurs if $X$ has two or more maximal
stabilizers, by the Hyperbolic Segment Condition
(\ref{hypsegment}). Otherwise, if $X$ has one maximal stabilizer and
there are no hyperbolic elements, then $G$ has a global fixed point by
Proposition \ref{titslemma}. In this case we define $X_G$ to be the tree
of fixed points of $G$, in keeping with the definition of $X_V$ above. 

\begin{lemma} \label{irrelevantedge}
If $e \in E(X) - E(X_G)$ separates $\partial_0 e$ from $X_G$, then $G_e =
G_{\partial_0 e}$. 
\end{lemma} 

\begin{proof}
Note that every $\gamma \in G$ maps $X_G$ to itself. Then if $\gamma$
fixes $\partial_0 e$, it also fixes the path from $\partial_0 e$ to
$X_G$, which contains $e$. Thus $G_{\partial_0 e} \subseteq G_e$. 
\end{proof}

\begin{definition} \label{telescope} 
A segment $(e_1, \ldots, e_n)$ is \emph{monotone} if either $G_{e_i} =
G_{\partial_0 e_i}$ for all $i$, or $G_{\overline{e}_i} = G_{\partial_0
\overline{e}_i}$ for all $i$. We say that $\partial_0 e_1$ is
\emph{minimal} and $\partial_1 e_n$ is \emph{maximal} in the first case,
and the reverse in the second case. An endpoint may be both minimal and
maximal; this occurs if and only if all stabilizers along the segment are
equal. We also declare every segment of length zero to be monotone. 

A segment $[v,w]$ is \emph{injective} if its interior maps injectively to
$G \backslash X$.  A segment $[v,w]$ is called a \emph{telescope} if
there is an edge $e\in [v,w]$ which separates $[v,w]$ into two monotone
segments, in which $\partial_0 e$ and $\partial_1 e$ are minimal. Note
that telescopes have positive length. 

Such an edge $e$ is called a $(v,w)$--minimal edge. Note that
$\overline{e}$ is also $(v,w)$--minimal. If $[v,w]$ is a
telescope then an edge $e\in [v,w]$ is $(v,w)$--minimal if and only if
$G_e = G_v \cap G_w$. The set of $(v,w)$--minimal edges forms a
nontrivial sub-segment of $[v,w]$, possibly equal to all of $[v,w]$. If
$[v,w] = [V,W]$ with $V,W\in M(X)$, a $(v,w)$--minimal edge is also called
a $(V,W)$--minimal edge. 

A \emph{proper telescope} is a telescope $[v,w]$ that is not monotone:
$G_v \not\subseteq G_w$ and $G_w \not\subseteq G_v$. 
\end{definition}

In the next definition we consider a pair of $G$--trees. Suppose that $Z$
is a proper $G$--tree and $X$ is obtained from $Z$ by an elementary
deformation. In the transition from $Z$ to $X$, various edges may have
stretched into telescopes. Also, vertices may have expanded into subtrees
whose vertices and edges have a common stabilizer (these are the
characteristic subtrees $X_U$). Other changes are likely to have occurred
as well. We would like to know that the pattern of these telescopes and
characteristic subtrees has the same combinatorics as $Z$. For this
reason we must eventually make the assumption that $Z$ is strongly
slide-free. 

The definition encodes what is meant by the pattern of telescopes having
the same combinatorics as $Z$. The somewhat elaborate requirements, taken
together, are invariant under elementary moves if $Z$ is strongly
slide-free. This is the content of Proposition \ref{telescopeprop}. The
requirements also imply that $Z$ can be recovered simply by performing
collapse moves. In this way Theorem \ref{mainrigidthm} is proved. 

\begin{definition} \label{telescoping} 
Let $X$ be a $G$--tree that possesses a maximal stabilizer. 
Suppose that $Z$ is a $G$--tree with vertex set $M(X)$, with the same
$G$--action (by conjugation). Note that $Z$ must be proper. We say that
$X$ is a \emph{proper telescoping} of $Z$ if: 
\begin{enumerate}
\renewcommand{\theenumi}{\roman{enumi}} 
\renewcommand{\labelenumi}{(\theenumi)} 
\item \label{e1} for every pair $V,W\in M(X)$ that bounds an edge in $Z$,
the segment $[V,W]$ is a telescope (call $[V,W]$ a
\emph{$Z$--telescope});
\item \label{e2} each $Z$--telescope $[V,W]$ is the union of two injective
telescopes $[V,w]$ and $[v,W]$, whose intersection $[v,w]$ is a proper
telescope, 
such that either $[v,w] = [V,W]$ or $[w,W] = t([v,V])$ for some $t \in
G$; 
\item \label{e3} for every $Z$--telescope $[V,W]$, the sub-telescope
$[v,w]$ given by (\ref{e2}) contains a $(V,W)$--minimal edge $e_{VW}$
which is not contained in any other $Z$--telescope (call $e_{VW}$ a
\emph{$Z$--minimal edge}).
\end{enumerate}
If we drop the requirement in (\ref{e2}) that $[v,w]$ is proper, then $X$
is called a \emph{telescoping} of $Z$. 
\end{definition}

The meaning of condition (\ref{e2}) is that each $Z$--telescope has one
of two configurations. Either it is proper and injective (the case $[v,w]
= [V,W]$), or its image in the quotient graph is a ``lasso,'' a path
which travels along a segment, around a circuit, then back along the
initial segment. Notice that the sub-telescope
$[v,w] \subseteq [V,W]$ is uniquely determined, as the
maximal injective telescope which is positioned symmetrically along
$[V,W]$. It maps to the ``circuit'' part of the lasso. Note that the
$Z$--minimal edge $e_{VW}$ in (\ref{e3}) is 
simultaneously a minimal edge for the telescopes $[V,W]$, $[V,w]$,
$[v,W]$, and $[v,w]$.

Note that every proper $G$--tree is a proper telescoping of itself. Also,
in any proper tree, telescopes have length $1$. Hence 
if $X$ is a telescoping of $Z$ and $X$ is proper then $X \cong Z$.

\begin{lemma} \label{minsubtree2} 
If $X$ is a telescoping of $Z$ then 
\begin{equation*} 
{\textstyle X_G \ = \ \bigcup_{U \in M(X)} \, (X_G \cap X_U) \ \cup \ 
\bigcup \, \{ Z\text{--telescopes} \} \, .}
\end{equation*}
\end{lemma}

\begin{proof}
First note that 
\begin{equation*} 
{\textstyle X_G \ \subseteq \ \bigcup_{U \in M(X)} \, X_U \ \cup \ 
\bigcup \, \{ Z\text{--telescopes} \} }
\end{equation*}
because the right hand side is a $G$--invariant subtree of $X$ (equal to
the subtree spanned by the $X_U$'s). To prove
the lemma it remains to show that every $Z$--telescope is contained in
$X_G$. Note that every $Z$--telescope has the property that its endpoint
stabilizers are strictly larger than any of its edge stabilizers (recall that
$[V,W]$ meets $X_V$ and $X_W$ only in its endpoints). The
Hyperbolic Segment Condition (\ref{hypsegment}) shows that such a segment is
contained in the axis of a hyperbolic element. Since $X_G$ contains every
hyperbolic axis, we are done. 
\end{proof}

\begin{lemma} \label{endpoints} 
Let $X$ be a telescoping of $Z$, and let $[V,W]$ and $[V',W']$
be $Z$--telescopes. If $\{V,W\} \cap \{V', W'\} = \emptyset$ then $[V,W]
\cap [V',W'] = \emptyset$. 
\end{lemma}

This implies that if two $Z$--telescopes meet then they have a common
endpoint. The lemma is a topological consequence of property
\ref{telescoping}(\ref{e3}).

\begin{proof}
Suppose that $\{V,W\} \cap \{V', W'\} = \emptyset$. 
The geometric edges $(V,W)$ and $(V',W')$
in $Z$ are separated by at least one edge. Consider the smallest segment
in $Z$ containing these two edges. Without loss of generality assume that
its vertices, in order, are $(W, V, U_1, \ldots, U_k, V', W')$ for some
$k \geq 0$. 

Let $T_i \subseteq X$ be the smallest tree containing 
\begin{equation} \label{chartrees} 
(X_W \cup X_V \cup X_{V'} \cup X_{W'}) \ \cup \ (X_{U_1} \cup \cdots \cup
X_{U_i}) 
\end{equation}
for each $i\geq 1$, and let $T_0$ be the span of $(X_W \cup X_V \cup
X_{V'} \cup X_{W'})$. Then each $T_i$ is the union of the characteristic
subtrees in \eqref{chartrees} and the $Z$--telescopes joining them. 
More specifically, for each $i\geq 1$ we have that 
\[T_{i+1} \ = \ T_i \ \cup \ [U_i, U_{i+1}] \ \cup \ X_{U_{i+1}},\]
as the right hand side is certainly a
(connected) tree containing the required subtrees. Similarly, $T_1 \ = \
T_0 \ \cup \ [V,U_1] \ \cup \ X_{U_1}$. 

Now suppose that $[V,W] \cap [V',W'] \not= \emptyset$. This implies that 
\[T_0 \ = \ (X_W \cup X_V \cup X_{V'} \cup X_{W'}) \ \cup \ [W,V] \ \cup
\ [V',W'].\] 
Together with the previous remarks we now have that 
\begin{multline} \label{chartree} 
T_k \ = \ (X_W \cup X_V \cup X_{V'} \cup X_{W'}) \ \cup \ (X_{U_1} \cup
\cdots \cup X_{U_k}) \\ \cup \ ([W,V] \cup [V',W']) \ \cup \ ([V,U_1]
\cup [U_1, U_2] \cup \cdots \cup [U_{k-1}, U_k]). 
\end{multline}
Consider the $Z$--telescope $[U_k,V']$. It is contained in $T_k$ and its
interior meets no $X_U$, $U \in M(X)$. Therefore equation
\eqref{chartree} implies that 
\begin{equation}
[U_k,V'] \ \subseteq \ ([W,V] \cup [V',W']) \ \cup \ ([V,U_1]
\cup [U_1, U_2] \cup \cdots \cup [U_{k-1}, U_k]).\tag*{} 
\end{equation}
Thus every edge of $[U_k,V']$ is contained in another $Z$--telescope,
contradicting property \ref{telescoping}(\ref{e3}). 
\end{proof}

\begin{lemma} \label{ssfproper} 
Let $X$ be a telescoping of $Z$, with $Z$ strongly slide-free. Then $X$
is a proper telescoping of $Z$. 
\end{lemma}

\begin{proof}
Let $[V,W]$ be a $Z$--telescope with injective sub-telescope $[v,w]$. If
$[v,w] = [V,W]$ then it is proper because $V$ and $W$ are maximal. If
$t([V,v]) = [W,w]$ for some $t\in G$ and $G_v \subseteq G_w$, then $G_v =
V \cap W$ because $[v,w]$ contains a $(V,W)$--minimal edge. As $t(v) = w$,
we also have $G_w = (V^t \cap W^t)$. Let $e\in E(Z)$ be the edge with
initial vertex $V$ and terminal vertex $W$. Then $te$ has initial vertex
$W$, and hence $\overline{e}, te \in E_0(W)$. Note that $G_{\overline{e}}
= V\cap W = G_v$ and $G_{te} = V^t \cap W^t = G_w$, so we have
$G_{\overline{e}} \subseteq G_{te}$. As $Z$ is strongly slide-free,
$\overline{e} \in G te = G e$. But then $e$ is inverted by an element of
$G$, a contradiction. Thus $G_v \subseteq G_w$ cannot hold, and similarly
$G_w \not\subseteq G_v$, so $[v,w]$ is proper.
\end{proof}

\begin{lemma} \label{zminedge2} 
Let $X$ be a telescoping of $Z$, with $Z$ strongly slide-free. Suppose
that $e$ is a $Z$--minimal edge with $G_e = G_{\partial_0 e}$. Then
$\partial_0 e$ is contained in exactly one $Z$--telescope.
\end{lemma}

\begin{proof}
Suppose that $e$ is $Z$--minimal for the $Z$--telescope $[V,W]$.  If
$\partial_0 e$ is contained in another $Z$--telescope, then by Lemma
\ref{endpoints} the two telescopes have a common endpoint. Suppose the
second $Z$--telescope is $[V,W']$. Then since $G_{\partial_0 e} =
G_e = V \cap W$ and $\partial_0 e \in [V,W']$, we must have $(V \cap W)
\supseteq (V \cap W')$, as $[V,W']$ is a telescope. 

Now since $Z$ is strongly slide-free, $V = V^{\delta}$ and $W' =
W^{\delta}$ for some $\delta \in G$. Then $\delta([V,W]) = [V,W']$,
fixing the endpoint in $X_V$. Since $\partial_0 e \in [V,W] \cap [V,W']$,
$\partial_0 e$ is fixed by $\delta$. But $\delta e \not= e$ as $e \not\in
[V,W']$, and hence $\delta \in G_{\partial_0 e} - G_{e}$, a
contradiction.
\end{proof}

The rest of this section deals with invariance of properties under
an elementary move from $X$ to $Y$. Without ambiguity, we will sometimes
use the same symbol to denote an edge or vertex of $X$ and its image in
$Y$, except for those edges and vertices involved in the move. 

\begin{lemma} \label{monotonecoll}
Let $(e_1, \ldots, e_n)$ be a monotone segment with $\partial_0 e_1$
minimal. Suppose that a collapse move is performed at $e\in E(X)$, where
$G_e = G_{\partial_0 e}$. Then exactly one of the following occurs:
\begin{enumerate}
\item \label{c1} the image of $(e_1, \ldots, e_n)$ is monotone with the
image of $\partial_0 e_1$ minimal, \nolinebreak or
\item \label{c2} $G_e \subsetneq G_{\partial_1 e}$ and there is a $\gamma
\in G$ such that $\partial_0 \gamma e = \partial_0 e_i$ for some $i$, and
$\gamma e \not= e_i, \overline{e}_{i-1}$. 
\end{enumerate}
\end{lemma}

\begin{proof} 
Let $v_i = \partial_0 e_i$ for $i = 1, \ldots, n$ and $v_{n+1} =
\partial_1 e_n$. Let $w_i$ be the image of $v_i$ under the collapse
move. Then the image segment $[w_1, w_{n+1}]$ is monotone with $w_1$
minimal if and only if $G_{w_i} \subseteq G_{w_{i+1}}$ for $i = 1,
\ldots, n$. Of course, $w_i$ and $w_{i+1}$ may be the same vertex for
various values of $i$.

(\ref{c1})$\Rightarrow \neg$(\ref{c2}): Assumption (\ref{c1}) implies
that $G_{w_i} \subseteq G_{w_{i+1}}$ for each $i$. Now suppose that
an element $\gamma \in G$ as in (\ref{c2}) exists. There are two
cases. First, if $G_{w_{i+1}} = 
G_{v_{i+1}}$ then 
$G_{\partial_1 \gamma e} = G_{w_i} \subseteq G_{v_{i+1}}$. Recall the
elementary fact that if $[x,y]$ is any segment with $G_x \subseteq G_y$
then $G_x \subseteq G_e$ for every $e \in [x,y]$. Thus, as $\gamma e \in
[\partial_1 \gamma e,v_{i+1}]$, we have that $G_{\partial_1 \gamma e}
\subseteq G_{\gamma e}$, contradicting (\ref{c2}). The second case is
that $G_{w_{i+1}} \supsetneq G_{v_{i+1}}$, meaning that $v_{i+1} =
\partial_0 \delta e$ for some $\delta\in G$, so $G_{w_{i+1}} =
G_{\partial_1 \delta e}$. Notice that $\gamma e \in [\partial_1 \gamma e,
\partial_1 \delta e]$, so again we have that $G_{\partial_1 \gamma e}
\subseteq G_{\gamma e}$, contradicting (\ref{c2}). 

$\neg$(\ref{c2})$\Rightarrow$(\ref{c1}): If $G_e = G_{\partial_1 e}$ then
no vertex stabilizer changes during the collapse move. Thus $G_{w_i} =
G_{v_i}$ for all $i$ and so $G_{w_i} \subseteq G_{w_{i+1}}$ for each
$i$. Otherwise, if $G_e \subsetneq G_{\partial_1 e}$, then the assumption
$\neg$(\ref{c2}) implies that for every $\gamma \in G$ and every $i\leq
n$, either $\partial_0 \gamma e \not= v_i$ or $\gamma e = e_i$ or $\gamma e =
\overline{e}_{i-1}$. 

If $v_i \not\in G \partial_0 e$ then $G_{w_i} = G_{v_i}$, so $G_{w_i}
\subseteq G_{v_{i+1}} \subseteq G_{w_{i+1}}$. If $v_i = \partial_0 \gamma
e$ with $i \leq n$ then either $\gamma e = e_i$ or $\gamma e =
\overline{e}_{i-1}$. If $\gamma e = 
\overline{e}_{i-1}$ then $\partial_1 \gamma e = v_{i-1}$ and so $G_{w_i}
= G_{v_{i-1}} \subseteq G_{v_{i+1}} \subseteq G_{w_{i+1}}$. If $\gamma e
= e_i$ then $\partial_1 \gamma e = v_{i+1}$, so $G_{w_i} = G_{v_{i+1}}
\subseteq G_{w_{i+1}}$. Combining all cases we have that $G_{w_i} \subseteq
G_{w_{i+1}}$ for $i = 1, \ldots, n$, and hence the segment $[w_1,
w_{n+1}]$ is monotone with $w_1$ minimal. 
\end{proof}

\begin{lemma} \label{monotoneexp} 
Let $(e_1, \ldots, e_n)$ be a segment and suppose that a collapse move is
performed at $e\in E(X)$, where $G_e = G_{\partial_0 e}$. Suppose that
the image of $(e_1, \ldots, e_n)$ is monotone with the image of
$\partial_0 e_1$ minimal. 

If $\overline{e}_n \not\in G e$ or $G_e = G_{\partial_1 e}$ then
$(e_1, \ldots, e_n)$ is monotone with $\partial_0 e_1$
minimal. Otherwise, $(e_1, \ldots, e_{n-1})$ is monotone with $\partial_0
e_1$ minimal. 
\end{lemma}

This lemma describes preservation of monotonicity of segments by
expansion moves, though it is expressed in terms of a collapse move. 

\begin{proof}
Let $v_i = \partial_0 e_i$ for $i = 1, \ldots, n$ and $v_{n+1} =
\partial_1 e_n$. Let $w_i$ be the image of $v_i$ under the collapse
move. Then $G_{v_i} \subseteq G_{w_i}$ for each $i$. 

For each $i\leq n$, $e_i$ and $e_{i+1}$ are not both in $G e$; this would
imply that $v_{i+1}\in (G \partial_0 e \cap G \partial_1 e)$, but this
set is empty because $e$ admits a collapse move. Similarly,
$\overline{e}_i$ and $\overline{e}_{i+1}$ are not both in $G e$. If
$\overline{e}_i, e_{i+1}\in G e$ then there is a $t\in G$ which fixes
$\partial_0 e_{i+1}$ but takes $e_{i+1}$ to $\overline{e}_i$. But then $t
\in G_{\partial_0 e_{i+1}} - G_{e_{i+1}}$, contradicting the assumption
that $G_e = G_{\partial_0 e}$, and so $\overline{e}_i$ and $e_{i+1}$ are
not both in $G e$. Finally suppose that $e_i, \overline{e}_{i+1} \in G
e$, for some $i < n-1$. Notice that since $e_i$ collapses to $w_{i+1}$,
$G_{w_{i+1}} = G_{v_{i+1}}$. The 
previous cases imply that $e_{i+2}, \overline{e}_{i+2} \not\in G e$ and
so $e_{i+2}$ is not collapsed, and has the same stabilizer after the
move. Then monotonicity of $[w_1, w_{n+1}]$ implies that $G_{w_{i+1}} =
G_{e_{i+2}}$. Hence $G_{v_{i+1}} = G_{e_{i+2}}$. There is a $t\in G$
fixing $\partial_1 \overline{e}_{i+1}$ ( $= v_{i+1}$) and taking
$\overline{e}_{i+1}$ to $e_i$. Then $t\in G_{e_{i+2}}$, and therefore $t$
fixes $\overline{e}_{i+1}$, as this edge lies between $v_{i+1}$ and
$e_{i+2}$. This contradicts $t \overline{e}_{i+1} = e_i$, and hence $e_i$
and $\overline{e}_{i+1}$ are not both in $G e$.

These facts imply that if $e_i$ is collapsed (ie, $e_i \in (G e \cup G
\overline{e})$) then $e_{i+1}$ is not collapsed, unless $i = n-1$ and
$e_{n-1}, \overline{e}_n \in G e$. Now consider an edge $e_i$. If $e_i$
is not collapsed then $G_{w_i} \subseteq G_{e_i}$ by monotonicity, and so
$G_{v_i} \subseteq G_{w_i} \subseteq G_{e_i}$.

If $e_i$ is collapsed then either $e_i \in G e$ or $\overline{e}_i \in
Ge$. In the first case $G_{v_i} \subseteq G_{e_i}$ holds. So suppose that
$\overline{e}_i \in G e$. Then $G_{w_{i+1}} = G_{w_i} = G_{v_i}$, and
$G_{e_i} = G_{v_{i+1}} \supseteq G_{e_{i+1}}$. Then, assuming that $i<n$,
the edge $e_{i+1}$ is not collapsed; the exception $i=n-1$, $e_i,
\overline{e}_{i+1}\in G e$, is ruled out because $\overline{e}_i \in
Ge$. Thus after the collapse $e_{i+1}$ is present, with stabilizer
$G_{e_{i+1}}$, and $G_{w_{i+1}} \subseteq G_{e_{i+1}}$ by
monotonicity. Hence $G_{v_i} \subseteq G_{e_{i+1}} \subseteq G_{e_i}$. 

Thus we have shown that $G_{v_i} \subseteq G_{e_i}$ for all $i<n$, and
for $i=n$ if $\overline{e}_n \not\in G e$. If $\overline{e}_n \in Ge$ and
$G_e = G_{\partial_1 e}$ then $G_{v_n} = G_{e_n}$ by equivariance. 
\end{proof}

\begin{lemma} \label{collapsible} 
Let $X$ be a telescoping of $Z$, with $Z$ strongly slide-free. Let
$[V,W]\subseteq X$ be a $Z$--telescope. Suppose that $e$ is a collapsible
edge of $X$, with $G e \cap [V,W] = \emptyset$. Then one of the following
occurs: 
\begin{enumerate}
\item \label{q1} $G \partial_0 e \cap [V,W] = \emptyset$, or 
\item \label{q2} $G \partial_1 e \cap [V,W] = \emptyset$, or 
\item \label{q3} $G \partial_0 e$ and $G \partial_1 e$ each meet $[V,W]$
in one endpoint, and the corresponding elements of $G e$ are in $X_V$ and
$X_W$. 
\end{enumerate}
\end{lemma}

\begin{proof}
Recall that collapsibility of $e$ means that $G_e = G_{\partial_0 e}$ and
$G \partial_0 e \cap G \partial_1 e = \emptyset$. Suppose that (\ref{q1})
and (\ref{q2}) do not hold, so that $[V,W]$ meets both $G \partial_0 e$
and $G \partial_1 e$. Relabeling if necessary, we have $\partial_0 e \in
[V,W]$ and $\partial_1 \gamma e \in [V,W]$ for some $\gamma \in G$. Note
that by Lemma \ref{minsubtree2} the endpoints of $e$ are in $X_G$, and
hence $e\in X_G$. For every $U \not= V,W$ in
$M(X)$, $X_U$ does not meet $[V,W]$ and therefore $e \not\in X_U$. If
$e\in X_V$ then $\partial_0 e$ is the endpoint $[V,W] \cap X_V$. The edge
$\gamma e$ is in $X_{(V^{\gamma})}$, and since $\partial_1 \gamma e \in
[V,W]$ we must have $X_{(V^{\gamma})} = X_W$. If $G \partial_0 e$ and $G
\partial_1 e$ do not meet the interior of $[V,W]$ then condition
(\ref{q3}) now holds. The case $e \in X_W$ is similar. 

Otherwise, we can assume that $e, \gamma e \not\in X_U$ for all $U \in
M(X)$. Then Lemma \ref{minsubtree2} implies that $e$ is in another
$Z$--telescope, which has a common endpoint with $[V,W]$ by Lemma
\ref{endpoints}. 

There are now two cases. Suppose first that the $Z$--telescope containing
$e$ has the form $[V,W']$. Then $\gamma e$ is in the $Z$--telescope
$\gamma ([V,W'])$, so Lemma \ref{endpoints} implies that $\{V,W\} \cap
\{V^{\gamma}, (W')^{\gamma}\} \not= \emptyset$.  Note that $e$ separates
$X_V$ from $\partial_1 e$, and so $\gamma e$ separates $X_{(V^{\gamma})}$
from $\partial_1 \gamma e \in [V,W]$. Hence $V^{\gamma} \not= V,W$ and so
$\{V,W\} \cap \{V^{\gamma}, (W')^{\gamma}\} = \{(W')^{\gamma}\}$. Thus
$(W')^{\gamma} = V$ or $(W')^{\gamma} = W$. 

The $Z$--telescope $[V,W']$ is separated by $e$ into two segments,
$[V,\partial_0 e]$ and $[\partial_1 e,W']$. The first is contained in the
$Z$--telescope $[V,W]$. The second is mapped by $\gamma$ either to
$[\partial_1 \gamma e,V] \subseteq [V,W]$ or to $[\partial_1 \gamma e,W]
\subseteq [V,W]$, and hence $[\partial_1 e,W']$ is contained in the
$Z$--telescope $\gamma^{-1}([V,W])$. 
Therefore $\{e,\overline{e}\}$ is the only geometric edge of $[V,W']$
that is not contained in another $Z$--telescope. It follows, from property
\ref{telescoping}(\ref{e3}), that $e$ is a $Z$--minimal edge for $[V,W']$. 
But the vertex $\partial_0 e$ is in two $Z$--telescopes, contradicting
Lemma \ref{zminedge2}.

The second case, when the $Z$--telescope containing $e$ has the form
$[V',W]$, yields a similar contradiction. 
\end{proof}

\begin{proposition} \label{telescopeprop} 
Let $Z$ be a strongly slide-free $G$--tree, and suppose that $X$ is a 
telescoping of $Z$. Let $Y$ be obtained from $X$ by an elementary
move. Then $Y$ is a telescoping of $Z$.
\end{proposition}

\begin{proof} 
First note that elementary moves do not change the set of maximal
stabilizers, and that $Z$ possesses maximal stabilizers, because it is
proper. Then $M(X) = M(Z)$ and $M(Y) = M(X)$. To avoid confusion, we
will decorate segments with subscripts to indicate which tree they are
in. 

\textbf{Step 1A}\qua {\slshape Invariance of
\textup{\ref{telescoping}(\ref{e1})} and
\textup{\ref{telescoping}(\ref{e3})} under collapses}\qua  Let $[V,W]_X$ be
a $Z$--telescope. Suppose the edge $e$ is collapsed, where $G_e =
G_{\partial_0 e}$. First we arrange that the $Z$--minimal edge $e_{VW}$
is not contained in $G e \cup G \overline{e}$. 

If $e_{VW} \in (Ge \cup G \overline{e})$, then we have $e_{VW} = \gamma
e$ for some $\gamma \in G$ (replacing $e_{VW}$ by its inverse if
necessary). Then $G_{\partial_0 e_{VW}} = G_{e_{VW}} = V\cap W$, and so
$\partial_0 e_{VW}$ is not an endpoint of $[V,W]_X$. Let $e'$ be the
unique edge in $(E_0(\partial_0 e_{VW}) \cap [V,W]_X) -
\{e_{VW}\}$. Since $[V,W]_X$ is a telescope, $G_{e'} = G_{e_{VW}}$, and
so $e'$ is a $(V,W)$--minimal edge. If $e'\not\in [v,w]$ then $\partial_0
e_{VW}$ is an endpoint of $[v,w]$, which implies that this telescope is not
proper. However, $X$ is a proper telescoping of $Z$ by Lemma
\ref{ssfproper}, so we conclude that $e'\in [v,w]$. Note also that $e'$
is not contained in any other $Z$--telescope of $X$, for then $\partial_0
e_{VW}$ would be in two $Z$--telescopes, contradicting Lemma
\ref{zminedge2}. Finally we note that $e' \not\in (G e \cup G
\overline{e})$. If $e' \in G \overline{e}$ then $\partial_0 e_{VW} \in G
\partial_0 e \cap G \partial_1 e$, but this set is empty as $e$ is
collapsible. If $e' \in G e$ then $\delta e' = e_{VW}$ for some $\delta
\in G$, but then $\delta \in (G_{\partial_0 e_{VW}} - G_{e_{VW}}) =
\emptyset$.

Thus, replacing $e_{VW}$ with $e'$ if necessary, there is a $Z$--minimal
edge $e_{VW} \not\in (Ge \cup G \overline{e})$. Then $e_{VW}$ is present
after the collapse move. We can orient $e_{VW}$ so that $\partial_0
e_{VW}$ separates $X_V$ from $e_{VW}$. Applying Lemma \ref{monotonecoll}
to the monotone segments $[V, \partial_0 e_{VW}]_X$ and $[\partial_1
e_{VW}, W]_X$, we find that $[V,W]_Y$ is a telescope with
$(V,W)$--minimal edge $e_{VW}$, unless $G_e \subsetneq G_{\partial_1 e}$ 
and $\partial_0 \gamma e$ is in the interior of $[V,W]_X$ for some
$\gamma \in G$, and $\gamma e \not\in [V,W]_X$. 

If this occurs, then note that $[V,W]_X \subseteq X_G$ by Lemma
\ref{minsubtree2}, and so $\gamma e \in X_G$ by Lemma
\ref{irrelevantedge}. Also, $\gamma e \not\in X_U$ for all $U\in M(X)$,
since no $X_U$ meets the interior of $[V,W]_X$. Hence, by Lemma
\ref{minsubtree2}, $\gamma e$ is contained in a $Z$--telescope. Lemma
\ref{endpoints} implies that this telescope has a common endpoint with
$[V,W]_X$. Suppose the telescope is $[V,W']_X$. Note that the subsegment
$[\partial_1 \gamma e, W'] \subseteq [V,W']_X$ cannot contain a
$(V,W')$--minimal edge, because $G_{\gamma e} \subsetneq G_{\partial_1
\gamma e}$. The segment $[V,\partial_0 \gamma e]$ is the
intersection of $Z$--telescopes $[V,W]_X$ and $[V,W']_X$, and hence it does
not contain a $Z$--minimal edge. Thus the
$(V,W')$--minimal edge promised by property \ref{telescoping}(\ref{e3})
can only be $\gamma e$ or $\gamma \overline{e}$. Hence, $\gamma e$ is a
$Z$--minimal edge for $[V,W']_X$ (and so is $\gamma \overline{e}$). Now
$\partial_0 \gamma e$ is in two $Z$--telescopes, contradicting Lemma
\ref{zminedge2}. A similar contradiction ensues if the $Z$--telescope
containing $\gamma e$ has the form $[V',W]_X$. 

Therefore $[V,W]_Y$ is a telescope with $(V,W)$--minimal edge
$e_{VW}$. It is clear that $e_{VW}$ is not contained in any other
$Z$--telescope of $Y$, since every $[V',W']_Y$ is the image of
$[V',W']_X$, and this property holds in $X$. Lastly, we need to know that
$e_{VW}$ is contained in the injective sub-segment $[v',w']_Y \subseteq
[V,W]_Y$ given by property \ref{telescoping}(\ref{e2}) of $Y$. In Step 2A
we show that $Y$ indeed satisfies \ref{telescoping}(\ref{e2}), and that
$[v',w']_Y$ is the image of $[v,w]_X$. Therefore $e_{VW} \in [v',w']_Y$. 

\textbf{Step 1B}\qua {\slshape Invariance of
\textup{\ref{telescoping}(\ref{e1})} under expansions}\qua  Suppose that an
expansion move is performed, creating $e$ with $G_e = G_{\partial_0
e}$. By Lemma \ref{monotoneexp}, $[V,W]_Y$ is a telescope, 
unless the following conditions hold: $\gamma e \in [V,W]_Y$ for some
$\gamma \in G$ with $\partial_0 \gamma e$ equal to one of the endpoints
of $[V,W]_Y$, and $G_e \subsetneq G_{\partial_1 e}$. But then
$G_{\partial_0 \gamma e} = G_{\gamma e} \subsetneq G_{\partial_1 \gamma
e}$, which contradicts maximality of $V$ or $W$ (one of which is equal to
$G_{\partial_0 \gamma e}$). 

\textbf{Step 2A}\qua {\slshape Invariance of
\textup{\ref{telescoping}(\ref{e2})} under collapses}\qua  Let $[V,W]_X$ be
a $Z$--telescope, and suppose that the edge $e$ is collapsed, where $G_e
= G_{\partial_0 e}$. The image of $[V,W]_X$ in $Y$ is $[V,W]_Y$, which is
a $Z$--telescope by Step 1A. Let $v'$ and $w'$ be the images of
$v$ and $w$ respectively. If $t([V,v]_X) = [W,w]_X$ then also
$t([V,v']_Y) = [W,w']_Y$. Note that $[V,w']_Y$ and $[v',W]_Y$ are
telescopes, since they are nontrivial sub-segments of the telescope
$[V,W]_Y$. Now, to establish \ref{telescoping}(\ref{e2}) for $Y$, it
remains to show that $[V,w']_Y$ and $[v',W]_Y$ are injective (we are not
required to show that $[v',w']$ is proper).

If $[V,w']_Y$ is not injective then two of its interior vertices are in
the same orbit. The preimages in $[V,w]_X$ of these vertices must be
incident to opposite ends of the geometric edge orbit $G \{e,
\overline{e}\}$.  By relabeling if necessary, we can arrange that
$\partial_0 e$ and $\partial_1 \gamma e$ are in the interior of $[V,w]_X$
for some $\gamma \in G$. Injectivity of $[V,w]_X$ implies that if $\delta
e \in [V,w]_X$ for some $\delta\in G$, then the endpoints of $\delta e$
are among the vertices $\partial_0 e$, $\partial_1 \gamma e$, $w$, and
the endpoint of $[V,w]_X$ in $X_V$. But $\delta e$ does not join
$\partial_0 e$ to $\partial_1 \gamma e$ because these vertices have
separate images in $Y$. Also $\delta e$ does not join an endpoint of
$[V,w]_X$ to one of $\partial_0 e$, $\partial_1 \gamma e$, because the
images of these vertices are in the interior of $[V,w']_Y$. Therefore $G
e \cap [V,w]_X = \emptyset$. By property \ref{telescoping}(\ref{e2}) we
have $G e \cap [w,W]_X = \emptyset$, and so $G e \cap [V,W]_X =
\emptyset$.  We now obtain a contradiction by applying Lemma
\ref{collapsible}, since none of its three conclusions holds. Thus,
$[V,w']_Y$ is injective. Similarly, $[v',W]_Y$ is injective. 

\textbf{Step 2B}\qua {\slshape Invariance of
\textup{\ref{telescoping}(\ref{e2})} under expansions}\qua 
Suppose that an expansion move is performed at $u$, creating the
edge $e$ with $G_e = G_{\partial_0 e}$. Consider the intersection of
$[V,W]_Y$ with $G e \cup G \overline{e}$. First we observe that
no two adjacent edges of $[V,W]_Y$ are in $G e \cup G \overline{e}$. 
For if this occurs, then the two edges cannot be coherently oriented.  
Their common endpoint would be in $G \partial_0 e \cap G \partial_1
e$, but this intersection is empty because $e$ is collapsible. Also, if
the two edges were incoherently oriented, then they would be related by
an element of $G$ fixing their common endpoint. This element would
violate the telescope property of $[V,W]_Y$, which holds by Step 1B. 

Thus, no two geometric edges in $[V,W]_Y \cap (G e \cup G \overline{e})$
collapse to the same vertex of $X$. The telescope $[V,W]_Y$ is now
obtained from $[V,W]_X$ by replacing each of various vertices in $G u
\cap [V,W]_X$ by a single geometric edge in $G \{e, \overline{e}\}$
(vertices are omitted if the corresponding elements of $G \{e,
\overline{e}\}$ are not in $[V,W]_Y$). Note that property
\ref{telescoping}(\ref{e2}) implies that $[V,W]_X$ contains at most two
vertices of $G u$, and if there are two, then they are in $[V,v]$ and
$[W,w]$ respectively, and are related by $t$. 

Suppose that $\gamma u$ is in the interior of $[V,v]_X$, and suppose that
$\gamma u$ expands to $\gamma e\in [V,v]_Y$. Let $e_1, e_2 \in
E_0(\gamma u)$ be the two neighboring edges of $[V,v]_X$. Then $e_1$ and
$e_2$ are separated by $\gamma e$ in $Y$, and $te_1$ and $te_2$ are in
$[V,W]_Y$, separated by $t \gamma e$. Hence $t \gamma e \in
[w,W]_Y$. Similarly if $\gamma u$ is the endpoint $[V,W]_X \cap
X_V$ and $\gamma e \in [V,W]_Y$, then $t \gamma e \in
[V,W]_Y$. In each of these cases, $G u$ does not meet $[v,w]$ by
property \ref{telescoping}(\ref{e2}), and so $[v,w]$ is unchanged during
the move, and hence $[V,w]_Y$ and $[v,W]_Y$ are injective. Also,
the relation $t([V,v]_Y) = [w,W]_Y$ holds.

Next suppose that $v= \gamma u$ for some $\gamma \in G$. If $w\not\in G
u$ then $[v,w] = [V,W]_X$ and $[V,W]_Y$ is injective. Otherwise, $w = tv$
for some $t$, and $t([V,v]) = [W,w]$. Recall that
\ref{telescoping}(\ref{e2}) implies that $G u$ meets $[V,W]_X$ in $v$ and
$w$ only. Then $[V,W]_Y$ meets $G\{e, \overline{e}\}$ in zero, one, or
two geometric edges. If $[V,W]_Y \cap G\{e, \overline{e}\} = \emptyset$,
let $v'$ and $w'$ be the vertices of $[V,W]_Y$ that are incident to
$G\{e, \overline{e}\}$, and which map to $v$ and $w$ respectively under
the collapse move from $Y$ to $X$. Then $[V,w']_Y$ and $[v',W]_Y$ are
injective and $t([V,v']_Y) = [W,w']_Y$. If $[V,W]_Y$ meets $G e$ in one
edge $\gamma e$ which collapses to $v\in X$, let $v'$ be the endpoint of
$\gamma e$ closest to $Y_V$. Let $w' \in [V,W]_Y$ be the vertex
corresponding to $w\in [V,W]_X$. Then $t([V,v']_Y) = [W,w']_Y$, and
$[V,w']_Y$ and $[v',W]_Y$ are injective.

If $[V,W]_Y$ contains two edges of $G e$, then call the edges $\gamma e$
and $\delta e$, where $\gamma e$ collapses to $\gamma u=v$ and $\delta e$
collapses to $t \gamma u=w$. Let $v''$ and $w''$ be the endpoints of
$[\gamma e, \delta e]$. Then $[V,v'']_Y$ and $[W,w'']_Y$ map bijectively to
$[V,v]_X$ and $[W,w]_X$ respectively under the collapse from $Y$ to $X$. 
Note that $t([V,v'']_Y) = [W,w'']_Y$.

We have that $v'' = \partial_i \gamma e$ for $i=0 \text{ or } 1$. Then $t
v'' = w''$ implies $\partial_i t \gamma e = w''$. Now if
$\partial_{(1-i)} \delta e = w''$, then $w'' \in (G \partial_0 e \cap G
\partial_1 e)$, contradicting the fact that $e$ is collapsible. Hence
$\partial_i \delta e = w'' = \partial_i t \gamma e$. This implies that
$\delta \gamma^{-1} t^{-1}$ fixes $w''$, because $\delta \gamma^{-1}
t^{-1} (t \gamma e) = \delta e$. Since $[w'',W]_Y$ is monotone with $w''$
minimal, $\delta \gamma^{-1} t^{-1}$ also fixes all of
$[w'',W]_Y$. Therefore $\delta \gamma^{-1}([V,\gamma e]) = [W,\delta
e]$. Setting $v' = \partial_{(1-i)} \gamma e$ and $w' = \partial_{(1-i)}
\delta e$, we now have that the injective segment $[v',w']_Y \subseteq
[V,W]_Y$ is the sub-telescope required by property
\ref{telescoping}(\ref{e2}), where the element $\delta \gamma^{-1}$ has
taken the place of $t$.

Lastly suppose that $G u$ meets the interior of $[v,w]$. Then
\ref{telescoping}(\ref{e2}) implies that $G u$ meets $[V,W]_X$ in no other
point. Then $[V,w]_Y$ and $[v,W]_Y$ are injective and
so \ref{telescoping}(\ref{e2}) is satisfied by $Y$. 

\textbf{Step 3}\qua {\slshape Invariance of
\textup{\ref{telescoping}(\ref{e3})} under expansions}\qua  Suppose that an
expansion move is performed. Let $e_{VW}$ be a $Z$--minimal edge for
$[V,W]_X$. In Step 1B we showed that $e_{VW}$ is $(V,W)$--minimal for
$[V,W]_Y$. It is contained in the segment $[v',w']_Y \subseteq
[V,W]_Y$, according to the description of $[v',w']_Y$ obtained in Step
2B. Also $e_{VW}$ is not contained in another $Z$--telescope after the
move, and so it is $Z$--minimal.
\end{proof}

\begin{definition}
Let $X$ be a $G$--tree which reduces to $Z$. We say that $X$ is
\emph{uniquely reducible to $Z$} if every maximal sequence of collapse
moves performed on $X$ results in the $G$--tree $Z$. 
\end{definition}

The following result includes Theorem \ref{mainrigidthm} of the
Introduction. 

\begin{theorem} \label{localrigidity} 
Let $X$ be a strongly slide-free $G$--tree, and let $Y$ be obtained from
$X$ by an elementary deformation. If $Y$ is proper, or cocompact and
reduced, then there is a unique isomorphism of $G$--trees $X \to
Y$. Otherwise $Y$ is a proper telescoping of $X$. If $Y$ is cocompact
then it is uniquely reducible to $X$.
\end{theorem}

\begin{proof}
Proposition \ref{telescopeprop} and Lemma \ref{ssfproper} together imply
that the property of being a proper telescoping of $Z$, when $Z$ is
strongly slide-free, is invariant under elementary deformations. Then
since $X$ is a proper telescoping of itself, $Y$ is a proper telescoping
of $X$. If $Y$ is proper then $Y\cong X$ as remarked earlier. Uniqueness
of the isomorphism follows from the fact that the only equivariant self
map of a proper $G$--tree is the identity. To prove this fact, let $\phi
\co Z \to Z$ be equivariant. If $Z$ is proper then the stabilizer of a
vertex does not fix any other vertex. Since equivariance requires
$G_{\phi(z)} \subseteq G_z$ for all $z$, the map is the identity. 

For the last conclusion of the theorem, let $\widehat{Y}$ be obtained
from $Y$ by a maximal sequence of collapse moves. Then $\widehat{Y}$ is
related to $X$ by a deformation, so it is a proper telescoping of $X$ by
\ref{telescopeprop} and \ref{ssfproper}. Since $\widehat{Y}$ is reduced,
every $X$--telescope in $\widehat{Y}$ has length $1$, by properties
\ref{telescoping}(\ref{e1}) and \ref{telescoping}(\ref{e2}). Thus
$\widehat{Y}\cong X$.  The case when $Y$ is cocompact and reduced follows
from the statement just proved. 
\end{proof}

\section{Rigidity}

In this section we combine Theorems \ref{mainthm} and \ref{localrigidity}
to generalize the rigidity theorem of Bass and Lubotzky given in
\cite{bass:rigidity}. This result is Corollary \ref{rigidcor}. The graph
of groups formulation of this result 
then yields Corollary \ref{uniquefactorcor} of the
Introduction.  We also apply the theory to generalized Baumslag--Solitar
trees (defined below), in Corollary \ref{gbscor}. 

First we have the following immediate consequence of Theorems
\ref{mainthm} and \ref{localrigidity}. It includes Corollary
\ref{rigiditycor} of the Introduction. 

\begin{theorem} \label{rigidity} 
Let $G$ be a group. Let $X$ and $Y$ be cocompact $G$--trees with the same
elliptic subgroups. Suppose that $X$ is strongly slide-free. Then $Y$ is
a proper telescoping of $X$, and hence is uniquely reducible to $X$. If
in addition $Y$ is reduced, then there is a unique isomorphism of
$G$--trees $X \to Y$. \endproof
\end{theorem}

If one imposes fixed point properties on the vertex stabilizers of a
$G$--tree, then the set of elliptic subgroups may be uniquely
determined. Or, this conclusion may hold for other reasons. The rest of
the results in this section concern such phenomena.

\begin{proposition} \label{properelliptic} 
Let $X$ and $Y$ be proper ${G}$--trees which define the same
partition of ${G}$ into elliptic and hyperbolic elements. Suppose in
addition that $X$ is locally finite. Then every vertex stabilizer of $X$
has a fixed point in $Y$. 
\end{proposition} 

\begin{proof} 
Call a subgroup of ${G}$ \emph{weakly elliptic} if its elements are all
elliptic. The vertex stabilizers of $X$ and $Y$ are maximal weakly
elliptic subgroups, by properness (and the Hyperbolic Segment
Condition). 

Let $H = {G}_x$ be a vertex stabilizer of $X$ and consider its action on
$Y$. By Proposition \ref{titslemma}, $H$ fixes either a vertex or an end
of $Y$. Suppose $H$ fixes $\varepsilon \in \partial Y$. Let $\Es$ be the
set of elliptic elements of ${G}$; then $H \subseteq ({G}_{\varepsilon}
\cap \Es)$. Note that ${G}_{\varepsilon}\cap \Es$ is a subgroup (and is
weakly elliptic), and hence $H = {G}_{\varepsilon} \cap \Es$ by
maximality. 

Suppose first that $\varepsilon$ is \emph{rational}, meaning that
${G}_{\varepsilon}$ contains a hyperbolic element $\gamma$
(ie, $\varepsilon$ is an endpoint of some hyperbolic axis). 
Then $\gamma H \gamma^{-1}$ is maximal weakly elliptic and is also
contained in ${G}_{\varepsilon} \cap \Es$, hence is equal to
$H$. Since $\gamma$ is hyperbolic, the subgroups $H$ and $\gamma H
\gamma^{-1}$ are the stabilizers of different vertices of $X$,
contradicting properness. 

Next assume that $\varepsilon$ is \emph{irrational}, so that
${G}_{\varepsilon} \subseteq \Es$, and hence $H =
{G}_{\varepsilon}$. Let $(y_i)_{i \geq 0}$ be any sequence of vertices
of $Y$ tending monotonically to $\varepsilon$. Then
${G}_{\varepsilon}$ is the increasing union $\bigcup_{i \geq 0}
\,({G}_{\varepsilon})_{y_i}$. We can assume that the inclusions
$({G}_{\varepsilon})_{y_i} \subseteq ({G}_{\varepsilon})_{y_{i+1}}$
are strict infinitely often, for otherwise ${G}_{\varepsilon}$ would
fix a vertex of $Y$. Each $({G}_{\varepsilon})_{y_i} =
({G}_{y_i})_{\varepsilon}$ is equal to the intersection $\bigcap_{j
\geq i} \,{G}_{y_j}$, and therefore  
\begin{equation} \label{end} 
{\textstyle {G}_{\varepsilon} \ = \ \bigcup_{i \geq 0} \, \bigcap_{j
\geq i} \, {G}_{y_j}.} 
\end{equation}\eject
Now we ask whether a vertex stabilizer of $X$ can have this structure.
Each stabilizer ${G}_{y_i}$ fixes either a vertex or an end of $X$, and
in the latter case the end must be irrational as above. Then, by
maximality of $G_{y_i}$, it must be equal to the stabilizer of that
vertex or end. Thus we have ${G}_{y_i} = {G}_{z_i}$ for some $z_i \in
V(X) \cup \partial X$. Note that $\{z_i\}$ is infinite, for otherwise
only finitely many of the inclusions
$({G}_{\varepsilon})_{y_i} \subseteq ({G}_{\varepsilon})_{y_{i+1}}$ would
be strict. 

Let $B_1(x)$ denote the ball of radius $1$ at $x$. As $X$ is locally
finite, $B_1(x)$ separates $X \cup \partial X$ into finitely 
many connected components, one of which meets $\{z_i\}$ in an infinite
set. Thus, by passing to a subsequence we can arrange that $\{z_i\}$ is
contained in a connected component of $(X \cup \partial X) -
B_1(x)$. Note that equation \eqref{end} remains valid after passing to any
subsequence. Now consider the subtree spanned by $\{z_i\}$. $B_1(x)$ is
disjoint from this subtree, so there is a vertex $x'$ separating
$x$ from every $z_i$. Then $G_x \cap G_{z_i} \subseteq G_{x'}$ for all
$i$, hence $G_x \cap \bigcap_{j\geq i} {G}_{z_j} \subseteq G_{x'}$ for
all $j$. But $G_x \cap \bigcap_{j\geq i} {G}_{z_j} = G_{\varepsilon} \cap
\bigcap_{j\geq i} {G}_{z_j} = \bigcap_{j\geq i}
{G}_{z_j} = \bigcap_{j\geq i} {G}_{y_j}$, so \eqref{end} yields 
${G}_{\varepsilon} \subseteq {G}_{x'}$. This contradicts properness of
$X$ as ${G}_{\varepsilon} = {G}_x$. 
\end{proof}

\begin{proposition} \label{Eproper} 
Let $X$ and $Y$ be cocompact $G$--trees which define the same partition of
$G$ into elliptic and hyperbolic elements. Suppose that $Y$ is proper,
and that every vertex stabilizer of $X$ has a fixed point in $Y$. Then
$X$ and $Y$ have the same elliptic subgroups. 
\end{proposition}

\begin{proof}
In fact, we show directly that $X$ and $Y$ are related by a deformation,
which implies the conclusion by Remark \ref{elemsubgroups}. The proof is
nearly the same as the proof of Theorem \ref{mainthm}. We can assume
without loss of generality that $X$ is 
reduced, by performing collapse moves. By Proposition \ref{morphism}
there exist a $G$--tree $X'$ obtained from $X$ by subdivision and an
equivariant morphism $\phi\co X' \to Y$. This morphism is surjective
because $Y$ is proper and hence minimal. Thus, by Proposition
\ref{bfprop}, it factors as $\phi = \psi \rho$ where $\rho \co X' \to Z$
is a finite composition of folds and $\psi \co Z \to Y$ induces an
isomorphism of quotient graphs. 

Since $Y$ is proper, Proposition \ref{preimage}(\ref{p1}) implies that
$\psi$ is an isomorphism. Thus, $\phi$ is a finite composition of folds. 
Each individual fold is equivariant and therefore it preserves ellipticity
of elements of $G$. Since their composition $\phi$ preserves
hyperbolicity, each fold must preserve hyperbolicity as well. Therefore
Proposition \ref{multifoldelem} applies to each fold and $X'$ and $Y$ are
related by an elementary deformation. Then $X$ and $Y$ are related by a
deformation. 
\end{proof}

\begin{corollary} \label{rigidcor} 
Let ${G}$ be a group. Let $X$ be a strongly slide-free ${G}$--tree and $Y$
a proper ${G}$--tree, both cocompact. Suppose that all vertex stabilizers
are unsplittable. If either
\begin{enumerate}
\item \label{s1} one of the trees has \textup{(FA)} vertex stabilizers, or
\item \label{s2} one of the trees is locally finite, 
\end{enumerate} 
then there is a unique isomorphism of ${G}$--trees $X \rightarrow Y$. 
\end{corollary}

\begin{proof}
The assumptions imply that $X$ and $Y$ define the same partition of $G$
into elliptic and hyperbolic elements. By Theorem \ref{rigidity} it
suffices to show that $X$ and $Y$ have the same elliptic subgroups. Note
that both trees are proper by Remark \ref{ssfremark}. Exchanging names if
necessary, assume that $X$ is the tree referred to in condition
(\ref{s1}) or (\ref{s2}). If (\ref{s1}) holds then $X$ and $Y$ have the
same elliptic subgroups by Proposition \ref{Eproper}. If (\ref{s2}) holds
then Propositions \ref{properelliptic} and \ref{Eproper} together yield
the same conclusion. 
\end{proof}

\begin{question}\label{question} 
Let $X$ and $Y$ be proper cocompact ${G}$--trees which define the same
partition of ${G}$ into elliptic and hyperbolic elements. Do $X$ and
$Y$ have the same elliptic subgroups? If so, then the hypotheses
(\ref{s1}) and (\ref{s2}) of Corollary \ref{rigidcor} could be dropped. 
\end{question} 

Corollary \ref{uniquefactorcor} of the Introduction is a special
case of the following result, which is a graph of groups interpretation
of Corollary \ref{rigidcor}. 

\begin{corollary} \label{uniquefactor} 
Let ${\bf A} = (A,\As,\alpha)$ and ${\bf B} = (B,\Bs,\beta)$ be graphs of
groups with finite underlying graphs and unsplittable vertex
groups. Suppose that ${\bf A}$ is strongly slide-free and ${\bf B}$
is proper, and that either: 
\begin{enumerate}
\item one of the graphs of groups has \textup{(FA)} vertex groups, or
\item one of the graphs of groups has finite index
edge-to-vertex inclusions. 
\end{enumerate}
Let $a_0$ and $b_0$ be basepoints in $A$ and $B$ respectively, and
suppose that there is an isomorphism $\psi \co \pi_1({\bf A},a_0) \to
\pi_1({\bf B},b_0)$. Then there exist an isomorphism $\Phi =
(\phi,(\gamma)) \co {\bf A} \to {\bf B}$ and an element $g \in \pi
[\phi(a_0),b_0] \subseteq \pi({\bf B})$ such that $\psi$ factors as 
$$\pi_1({\bf A},a_0) \stackrel{\Phi_{a_0}}{\longrightarrow} \pi_1({\bf
B},\phi(a_0)) \stackrel{\ad(g)}{\longrightarrow} \pi_1({\bf B},b_0).$$ 
In particular, ${\bf A}$ and ${\bf B}$ are isomorphic as graphs of
groups. 
\end{corollary}

Recall that our notion of isomorphism for graphs of groups is taken from
\cite[\S 2]{bass:covering}, and is more inclusive than the naive
definition (cf Corollary \ref{herrlichcor}). For the additional notation
and definitions we refer to \cite{bass:covering}. The conclusion and
proof of \ref{uniquefactor} are the same as those of Corollary 4.5 in
\cite{bass:rigidity}.

\begin{proof}
Let $G = \pi_1({\bf A},a_0)$, $X = (\widetilde{{\bf A},a_0})$, and $Y =
(\widetilde{{\bf B},b_0})$. Then $X$ and $Y$ are $G$--trees (via $\psi$
in the case of $Y$), and they satisfy the hypotheses of Corollary
\ref{rigidcor}. Hence there is an isomorphism of $G$--trees $\phi \co X
\to Y$. By \cite[4.2--4.5]{bass:covering} there is an isomorphism $\Phi =
(\phi,(\gamma)) \co {\bf A} \to {\bf B}$. The element $g\in
\pi[\phi(a_0),b_0]$ and factorization $\psi = \ad(g) \circ
\Phi_{a_0}$ exist by \cite[\S 3 and Theorem 4.1]{bassjiang}. 
\end{proof}

It turns out that a more restrictive variant of Corollary \ref{rigidcor}
can be proved without using the deformation theory of section
\ref{equivsection}. This version does not assume cocompactness. Following
\cite{bass:rigidity}, a ${G}$--tree is \emph{strict} if no edge
stabilizer contains the stabilizer of another geometric edge or vertex.
Strict ${G}$--trees are strongly slide-free, but the converse does not
hold. For example, a generalized Baumslag--Solitar tree that is not a
point (see definition \ref{gbs} below) cannot be strict, but it can
easily be slide-free or strongly slide-free.

\begin{theorem}
Let ${G}$ be a group and let $X$ and $Y$ be strict
${G}$--trees. If all vertex stabilizers have property \textup{(FA)} then there
is a unique isomorphism of ${G}$--trees $X \to Y$. 
\end{theorem}

\begin{proof}
Since both trees have (FA) vertex stabilizers, they have the same
elliptic subgroups. Then, as each tree is proper, the vertex stabilizers
in each tree are precisely the maximal elliptic subgroups. Thus the trees
have the same vertex stabilizers. By Proposition (3.5)(d) of
\cite{bass:rigidity} the trees are canonically isomorphic. It is here
where strictness of both trees is required. 
\end{proof}

The essential content of the proposition cited above is the following. 
If $X$ is a strict ${G}$--tree then it is proper, and hence 
vertices correspond bijectively with their stabilizers. Then, using 
strictness, one finds that two vertices bound an edge if and only if the
intersection of their stabilizers is maximal among pairwise intersections
of vertex stabilizers. Thus $X$ is intrinsically determined by its set of
vertex stabilizers. The property of strictness seems to be the most
general condition under which this conclusion holds. 

\begin{definition} \label{gbs} 
A \emph{generalized Baumslag--Solitar tree} is a $G$--tree in which all
edge and vertex stabilizers are infinite cyclic. The groups $G$ which
arise are called \emph{generalized Baumslag--Solitar groups}.  Examples
include the classical Baumslag--Solitar groups, torus knot groups, and
finite index subgroups of these groups. 
\end{definition}

Note that the preceding results do not apply to generalized
Baumslag--Solitar trees because ${\Z}$ is not unsplittable. The results
of \cite{bass:rigidity} also do not apply. However, the elliptic
subgroups of generalized Baumslag--Solitar groups are easily determined,
and we have the following result. 

\begin{corollary} \label{gbscor} 
Let $X$ and $Y$ be cocompact generalized Baumslag--Solitar trees with
group $G$. Assume that each tree has a minimal subtree not contained in a
line. Then there is an elementary deformation relating $X$ to $Y$. If $X$
is strongly slide-free, then $Y$ is uniquely reducible to $X$. If in
addition $Y$ is reduced, then there is a unique isomorphism $X \to Y$. 
\end{corollary}

The assumption on minimal subtrees is directly analogous to the Kleinian
groups property of being \emph{non-elementary}: that the limit set has
more than two points. If a $G$--tree $X$ has a minimal subtree, then the
boundary of this subtree is the limit set of $X$. In this way the two
properties are seen to coincide.  If a cocompact generalized
Baumslag--Solitar tree does not have this property then the group is
isomorphic to ${\Z}$, ${\Z} \times {\Z}$, or the Klein bottle group. 

\begin{proof}
First note that $G$ is torsion free: every hyperbolic element has
infinite order and every elliptic element is contained in an infinite
cyclic subgroup. Also note that for any $\gamma\in G$ and $i\not=0$, the
element $\gamma^i$ is hyperbolic if and only if $\gamma$ is, and when
this occurs they have the same axis.  

We will say that two subgroups $H, H' \subseteq G$ are
\emph{commensurable} if $H \cap H'$ 
has finite index in both $H$ and $H'$. Commensurability is an equivalence
relation on the subgroups of $G$. 
We say that two nontrivial
elements $\gamma, \delta\in G$ are commensurable if the infinite
cyclic subgroups they generate are commensurable. This occurs if and only
if $\gamma^i = \delta^j$ for some $i,j \not= 0$. Thus, commensurable
hyperbolic elements have the same axis. 

In a generalized Baumslag--Solitar tree every edge stabilizer has finite
index in its neighboring vertex stabilizers. This implies that all edge
and vertex stabilizers are in a single commensurability class. Because
these subgroups are infinite cyclic, this class contains all nontrivial
subgroups of vertex stabilizers. Therefore every nontrivial elliptic
element is commensurable with all of its conjugates. 

Now suppose that an element $\gamma \in G$ is commensurable with its
conjugates. If $\gamma$ is hyperbolic then let $L_{\gamma}$ be its
axis. For every $\delta \in G$, $\delta \gamma \delta^{-1}$ is
commensurable with $\gamma$ and hence its axis $L_{\delta \gamma
\delta^{-1}}$ is equal to $L_{\gamma}$. But $L_{\delta \gamma
\delta^{-1}} = \delta L_{\gamma}$, and so we have shown that $L_{\gamma}$
is $G$--invariant. Hence $L_{\gamma}$ contains any minimal subtree. 

These remarks imply that if a minimal subtree of a generalized
Baumslag--Solitar tree is not contained in a line, then the nontrivial
elliptic elements are exactly those elements which are commensurable with
all of their conjugates. Thus, $X$ and $Y$ define the same partition of
$G$ into elliptic and hyperbolic elements. Then, as all vertex
stabilizers have property (E) (by finite generation), Corollary
\ref{mainthmcor} implies that the two trees are related by a
deformation. The other conclusions follow from Theorems \ref{mainthm} and
\ref{rigidity}. 
\end{proof}

\section{Geometric Rigidity} \label{geomsec} 

In this section we examine the preceding results from the
coarse geometric point of view. We find that deformation-equivalence of
$G$--trees can be characterized geometrically, as stated in condition
(\ref{a3}) of Theorem \ref{maindefthm}. This result, combined with
Theorem \ref{mainrigidthm}, yields the quasi-isometric rigidity theorem
given in Corollary \ref{geomrigiditycor}. 

\begin{definition}
A map $f\co X\rightarrow Y$ between metric spaces $X$ and $Y$ is a
\emph{quasi-isometry} if there exist constants $K \geq 1$, $C\geq0$ such
that for every $x, x' \in X$, 
\[ K^{-1} d(x,x') - C \ \leq \ d(f(x),f(x')) \ \leq \ K d(x,x') + C, \] 
and if $d(y, f(X)) < C$ for every $y\in Y$. 

If $X$ and $Y$ admit ${G}$--actions by isometries 
then a map $f\co X \rightarrow Y$ is
\emph{$L$--equivariant} if for every $x\in X$ and $\gamma \in {G}$,
$d(f(\gamma x), \gamma f(x)) < L$. If $f$ is $L$--equivariant for
some $L$ then we say that $f$ is \emph{coarsely equivariant}. 

If $f \co X\rightarrow Y$ is a quasi-isometry then an 
\emph{$M$--quasi-inverse} is a quasi-isometry $g\co Y \rightarrow X$ such
that $d(x,g(f(x))) < M$ and $d(y,f(g(y))) < M$ for every $x\in X$ and
$y\in Y$. 
\end{definition}

\begin{lemma} \label{quasiinverse}
Let $f\co X\rightarrow Y$ be a $(K,C)$--quasi-isometry. Then there exists a
$(K, 3KC)$--quasi-isometry $g\co Y\rightarrow X$ which is a
$2KC$--quasi-inverse of $f$. 
If $X$ and $Y$ admit ${G}$--actions by isometries and $f$ is
$L$--equivariant, then any $M$--quasi-inverse of $f$ is $K(2M + C +
L)$--equivariant.
\end{lemma}

\begin{proof}
We define $g$ by choosing $g(y)$ to be any point of $X$ with 
$d(f(g(y)),y) < C$. It is straightforward to check that $g$ is a $(K,
3KC)$--quasi-isometry, and a $2KC$--quasi-inverse of $f$. 
Next let $g$ be any $M$--quasi-inverse to $f$, and take $\gamma \in {G}$
and $y\in Y$. Then 
\begin{multline*}
K^{-1} d(g(\gamma y), \gamma g(y)) - C \ \leq \ d(f(g( \gamma y)),f(\gamma
g(y))) \\
\leq \ d(f(g( \gamma y)), \gamma y) + d(\gamma y, \gamma f(g(y))) +
d(\gamma f(g(y)), f( \gamma g(y))).
\end{multline*}
The middle term on the second line is equal to $d(y, f(g(y)))$ and so
the second line is bounded by $2M+L$. 
Hence $d(g(\gamma y), \gamma g(y)) \leq K(2M + C + L)$ as desired. 
\end{proof}

Theorem \ref{mainthm} can now be extended as follows, completing the
proof of Theorem \ref{maindefthm}. 

\begin{theorem} \label{qicor}
Let ${G}$ be a group, and let $X$ and $Y$ be cocompact
${G}$--trees. The following conditions are equivalent. 
\begin{enumerate}
\item \label{z1} $X$ and $Y$ have the same elliptic subgroups. 
\item \label{z2} There exists an elementary deformation $X \to Y$. 
\item \label{z3} There exists a coarsely equivariant quasi-isometry $\phi
\co X \rightarrow Y$. 
\end{enumerate}
\end{theorem}

Note in particular the implication (\ref{z3})$\Rightarrow$(\ref{z2}). It
transforms a coarse, approximate relationship between ${G}$--trees
into a combinatorially precise one. We also remark that the implications
(\ref{z2})$\Rightarrow$(\ref{z3}) and (\ref{z3})$\Rightarrow$(\ref{z1}) 
do not require cocompactness. 

\begin{proof} The implication (\ref{z1})$\Rightarrow$(\ref{z2}) is given
by Theorem \ref{mainthm}. For (\ref{z2})$\Rightarrow$(\ref{z3}), we claim 
that a collapse move is itself an equivariant quasi-isometry. Then Lemma
\ref{quasiinverse} implies that the relation defined by condition
(\ref{z3}) is an equivalence relation. Since the relation of condition
(\ref{z2}) is generated by collapse moves, and the implication
(\ref{z2})$\Rightarrow$(\ref{z3}) holds for these moves, it must hold 
in general. 

To prove the claim, we show that a collapse move is a $(3,
2/3)$--quasi-isometry. The essential point is that if $e$ can be
collapsed, then the connected components of $G e$ have diameter at most
$2$. Then, if $(e_1, \ldots, e_d)$ is a geodesic of length $d$, then at
most $(2/3)(d+1)$ of its geometric edges are in $G
\{e,\overline{e}\}$. Thus the distance $d'$ between the endpoints after
the collapse satisfies
\[ (1/3)(d - 2) \ \leq \ d' \ \leq \ d.\]
Now consider a collapsible edge $e$ (with $G_e = G_{\partial_0 e}$) and a
translate $\gamma e$ which is adjacent to $e$. The endpoints of $e$ are
in different $G$--orbits, so either $\partial_0 e = \partial_0 \gamma e$
or $\partial_1 e = \partial_1 \gamma e$. In the former case we have that
$\gamma \in G_{\partial_0 e} - G_e$, a contradiction. Hence the latter
holds. It follows that every
edge in the connected component of $G
e$ containing $e$ is incident to $\partial_1 e$. Hence this component has
diameter at most $2$.

For (\ref{z3})$\Rightarrow$(\ref{z1}) we will establish the fact that if
$f\co X \rightarrow Y$ is any coarsely equivariant map of
${G}$--trees, then every vertex stabilizer of $X$ fixes a vertex of
$Y$. We then apply this fact both to $\phi$ and to a quasi-inverse
$\psi\co Y \rightarrow X$, which exists and is coarsely equivariant by
Lemma \ref{quasiinverse}. Then $X$ and $Y$ have the same elliptic
subgroups.

Suppose $f\co X \rightarrow Y$ is $L$--equivariant. If $\gamma \in {G}_x$
for some $x\in V(X)$, then $d(\phi(x), \gamma \phi(x)) = d(\phi(\gamma
x), \gamma \phi(x)) < L$. Thus $T$, the subtree spanned by ${G}_x
\phi(x)$, is a bounded subtree of $Y$. Now every element $\gamma$ of
${G}_x$ has a fixed point in $Y$, since if it had positive translation
length then the subset $\{\gamma^i \phi(x)\}_{i>0} \subseteq T$ would be
unbounded. Note also that a fixed point of $\gamma$ lies on the path
$[\phi(x), \gamma \phi(x)] \subseteq T$ and so every element of ${G}_x$
has a fixed point in $T$.

By Proposition \ref{titslemma} either ${G}_x$ has a fixed point in 
$Y$ or there is an end $\varepsilon \in \partial Y$ fixed by
${G}_x$. Suppose the latter holds. Let $y$ and $y'$ be vertices of
$T$. The subtree spanned by $y$, $y'$, and $\varepsilon$ meets $T$ in a
subtree containing $[y,y']$, so the intersection $[y,\varepsilon) \cap
[y', \varepsilon) \cap T$ is nonempty. As $T$ is bounded, there is a 
nearest vertex $z \in [y,\varepsilon) \cap [y', \varepsilon) \cap T$ to
$\varepsilon$, which is also the nearest vertex of $[y,\varepsilon) \cap
T$ to $\varepsilon$. Fixing $y$, this shows that $z\in [y', \varepsilon)$
for every $y' \in T$. Therefore ${G}_x$ fixes $z$ as each of its
elements fixes both $\varepsilon$ and some vertex of $T$. 
\end{proof}

\begin{proposition} \label{finitedist} 
Let ${G}$ be a group and $X$ a proper ${G}$--tree. Let $\phi \co V(X)
\to V(X)$ be an $L$--equivariant map. Then $\dist(\phi,\id) < L/2$. In
particular, if $\phi$ is equivariant then $\phi = \id$. 
\end{proposition}

\begin{proof}
Let $x$ be a vertex of $X$. We wish to show that $d(x, \phi(x)) <
L/2$. For each $\gamma \in {G}_x$ we have that 
$d(\phi(x), \gamma \phi(x)) = d(\phi(\gamma x), \gamma \phi(x)) < L$. Let
$m_{\gamma}$ be the vertex defined by 
\[m_{\gamma} \ = \ [x,\phi(x)] \cap [\phi(x), \gamma \phi(x)] \cap
[\gamma \phi(x), x]. \]
That is, $m_{\gamma}$ is the midpoint of the tripod spanned by $x$,
$\phi(x)$, and $\gamma \phi(x)$. Note that $\gamma$ fixes the
segment $[x,m_{\gamma}] \subseteq [x,\phi(x)]$. Suppose that $m_{\gamma}
\not= x$ for every $\gamma \in {G}_x$. Then ${G}_x$ fixes an
interior point of $[x, \phi(x)]$. However, as $X$ is proper, the
fixed point set of ${G}_x$ is $\{x\}$, and therefore $m_{\gamma} = x$
for some $\gamma \in {G_x}$. Then $x$ is the midpoint of
$[\phi(x), \gamma \phi(x)]$, a path of length less than $L$. 
\end{proof}

The following result is Corollary \ref{geomrigiditycor} of the
Introduction. 

\begin{theorem} \label{geomrigidity} 
Let $G$ be a group. Let $X$ and $Y$ be reduced cocompact $G$--trees, with
$X$ strongly slide-free. If $\phi \co X \to Y$ is a coarsely
equivariant quasi-isometry, then there is a unique equivariant isometry
from $X$ to $Y$, and it has finite distance from $\phi$. 
\end{theorem}

\begin{proof}
The first conclusion follows directly from Theorems \ref{qicor} and
\ref{localrigidity}. The second conclusion follows from
Proposition \ref{finitedist}. 
\end{proof}


\begin{thebibliography}

\bibitem{bass:covering}
\textbf{Hyman Bass}, \emph{Covering theory for graphs of groups}, J. Pure Appl.
  Algebra 89 (1993) 3--47

\bibitem{bassjiang}
\textbf{Hyman Bass}, \textbf{Renfang Jiang}, \emph{Automorphism groups of tree
  actions and of graphs of groups}, J. Pure Appl. Algebra 112 (1996) 109--155

\bibitem{bass:treelat}
\textbf{Hyman Bass}, \textbf{Ravi Kulkarni}, \emph{Uniform tree lattices}, J.
  Amer. Math. Soc. 3 (1990) 843--902

\bibitem{bass:rigidity}
\textbf{Hyman Bass}, \textbf{Alexander Lubotzky}, \emph{Rigidity of group
  actions on locally finite trees}, Proc. London Math. Soc. 69 (1994)
  541--575

\bibitem{bestvina:accessibility}
\textbf{Mladen Bestvina}, \textbf{Mark Feighn}, \emph{Bounding the complexity
  of simplicial group actions on trees}, Invent. Math. 103 (1991) 449--469

\bibitem{chiswell}
\textbf{I\,M Chiswell}, \emph{The {G}rushko--{N}eumann theorem}, Proc. London
  Math. Soc. 33 (1976) 385--400

\bibitem{cullermorgan}
\textbf{Marc Culler}, \textbf{John~W Morgan}, \emph{Group actions on {${\mathbb
  R}$}--trees}, Proc. London Math. Soc. 55 (1987) 571--604

\bibitem{dunwoody:folding}
\textbf{M\,J Dunwoody}, \emph{Folding sequences}, from: ``The {E}pstein
  birthday schrift'', Geometry and Topology Monographs, 1
  (1998)  139--158

\bibitem{forester:jsj}
\textbf{Max Forester}, \emph{On uniqueness of {JSJ} decompositions of finitely
  generated groups}, preprint (2001) {\tt arXiv:math.GR/0110176}

\bibitem{herrlich}
\textbf{Frank Herrlich}, \emph{Graphs of groups with isomorphic fundamental
  group}, Arch. Math. 51 (1988) 232--237

\bibitem{kurosh:groups}
\textbf{A\,G Kurosh}, \emph{The theory of groups}, Chelsea Publishing Co., New
  York (1960)

\bibitem{msw:announcement}
\textbf{Lee Mosher}, \textbf{Michah Sageev}, \textbf{Kevin Whyte},
  \emph{Quasi-actions on trees, Research announcement}, preprint (2000)
  {\tt arXiv:math.GR/0005210}

\bibitem{serre:trees}
\textbf{Jean-Pierre Serre}, \emph{Trees}, Springer-Verlag (1980)

\bibitem{stallings:foldings}
\textbf{John~R Stallings}, \emph{Foldings of ${G}$--trees}, from: ``Arboreal
  group theory (Berkeley, CA, 1988)'', Springer, New York (1991)  355--368

\bibitem{tits:treeauto}
\textbf{J Tits}, \emph{Sur le groupe des automorphismes d'un arbre}, from:
  ``Essays on Topology and Related Topics: Memoires d\'edi\'es \`a {G}eorges de
  {R}ham'', (A Haefliger, R Narasimhan, editors), Springer (1970)  188--211

\end{thebibliography}
\end{document}